\documentclass[preprint,3p,12pt,final]{elsarticle}

\graphicspath{{figure/}}

\usepackage[utf8]{inputenc}
\usepackage[ngerman,english]{babel}
\usepackage[T1]{fontenc}
\usepackage{marvosym} 

\usepackage{graphicx}
\usepackage{caption}
\usepackage{subcaption}
\usepackage{xcolor}
\usepackage{tabularx}

\setkeys{Gin}{draft=false}

\usepackage{mathtools}
\usepackage{amssymb}
\usepackage{amsfonts}
\usepackage{amsbsy} 
\usepackage{amsopn} 
\usepackage{amstext} 
\usepackage{amsxtra} 
\usepackage{bbm} 
\usepackage{dsfont} 

\usepackage{cancel} 
\usepackage{amscd} 
\usepackage[all]{xy} 
\usepackage{esvect} 
\usepackage{esint} 
\usepackage{empheq} 

\usepackage{eso-pic,picture}
\usepackage{lscape}
\usepackage{blindtext}
\usepackage{listings}
\usepackage{url}
\usepackage{textcomp} 
\usepackage{exscale} 
\usepackage{xspace} 
\usepackage{scalerel} 
\usepackage{stackengine} 
\usepackage{hhline}
\usepackage{lineno}
\usepackage{fancyhdr}
\usepackage{longtable}

\usepackage{hyphenat} 
\usepackage{testhyphens} 

\input ushyphex
\PassOptionsToPackage{svgnames,dvipsnames}{xcolor}
\usepackage{tcolorbox}
\tcbuselibrary{skins,breakable}

\usepackage{lmodern}
\usepackage{microtype}

\usepackage[numbers]{natbib}

\usepackage{amsthm} 

\newtheorem{theorem}{Theorem}

\addto\captionsenglish{\renewcommand{\appendixname}{Appendix }}

\usepackage[notref,notcite]{showkeys} 
\usepackage[mode=multiuser]{fixme} 
\fxusetheme{color}
\FXRegisterAuthor{hh}{envhh}{HH}
\FXRegisterAuthor{ar}{envar}{AR}
\FXRegisterAuthor{dk}{envdk}{DK}
\fxusetheme{color}

\numberwithin{equation}{section}

\RequirePackage[l2tabu]{nag}

\usepackage{varioref} 
\usepackage[colorlinks=true,citecolor=blue,
urlcolor=blue,linkcolor=blue]{hyperref} 
\usepackage[capitalise]{cleveref} 

\crefname{appendix}{}{}
\crefname{listing}{Algorithm}{Algorithms}
\crefname{equation}{}{}

\newcommand{\dt}{\ensuremath{\,\mathrm{d} t}}
\newcommand{\dx}{\ensuremath{\,\mathrm{d}\vec x}}

\renewcommand{\d}{\,\mathrm{d}} 

\newcommand{\pd}[2]{\frac{\partial #1}{\partial #2}}
\newcommand{\td}[2]{\frac{\d #1}{\d #2}}

\newcommand{\bx}{\vec{x}}
\newcommand{\bF}{\vec{f}}

\newcommand{\bc}{\vec{c}}

\newcommand{\beq}{\begin{equation}}
\newcommand{\eeq}{\end{equation}}

  \newdefinition{rmk}{Remark}

\newcommand{\N}{\ensuremath{\mathbb{N}}}

\newcommand{\R}{\ensuremath{\mathbb{R}}}

\renewcommand{\phi}{\ensuremath{\varphi}}

\renewcommand{\Phi}{\ensuremath{\varPhi}}

\renewcommand{\epsilon}{\ensuremath{\varepsilon}}

\renewcommand{\rho}{\ensuremath{\varrho}}

\newcommand{\timeMesh}{\ensuremath{\mathcal I_{\Delta t}}}

\newcommand{\interpol}{\ensuremath{\mathcal I}}

\renewcommand{\vec}[1]{\ensuremath{\mathbf{#1}}}

\newcommand{\ie}{, \mbox{i.\,e.,}~}

\newcommand{\wrt}{\mbox{w.~r.~t.~}}

\allowdisplaybreaks 





\newcommand{\doi}[1]{\textrm{\textsc{doi:}} \href{http://dx.doi.org/#1}{\nolinkurl{#1}}}

\journal{}

\begin{document}
\begin{frontmatter}

\title{Limiter-based entropy stabilization of semi-discrete
and fully discrete schemes for nonlinear hyperbolic problems}

\author[TUD]{Dmitri Kuzmin}
\ead{kuzmin@math.uni-dortmund.de}

\author[TUD]{Hennes Hajduk\corref{cor}}
\ead{hennes.hajduk@math.tu-dortmund.de}
\cortext[cor]{Corresponding author}

\address[TUD]{Institute of Applied Mathematics (LS III), TU Dortmund
University,\\ Vogelpothsweg 87, D-44227 Dortmund, Germany}

\author[LUT]{Andreas Rupp}
\ead{andreas@rupp.ink}

\address[LUT]{School of Engineering Science, Lappeenranta--Lahti University of Technology (LUT),\\ P.O.\ Box 20, FI-53851 Lappeenranta, Finland}

\begin{abstract}
  The algebraic flux correction (AFC) schemes presented in this
  work constrain a standard continuous finite element discretization
  of a nonlinear hyperbolic problem to satisfy relevant maximum
  principles and entropy stability conditions. The desired
  properties are enforced by applying a limiter to antidiffusive
  fluxes that represent the difference between the high-order
  baseline scheme and a property-preserving approximation of
  Lax--Friedrichs type. In the first step of the 
  limiting procedure, the given target fluxes are adjusted in
  a way that guarantees preservation of local and/or global
  bounds. In the second step,
  additional limiting is performed, if necessary,
  to ensure the validity of fully
  discrete and/or semi-discrete entropy inequalities. The limiter-based
  entropy fixes considered in this work are applicable to finite
  element discretizations of scalar hyperbolic equations and
  systems alike. The underlying inequality constraints
  are formulated using Tadmor's entropy stability theory.

  The proposed limiters impose
  entropy-conservative or entropy-dissipative
  bounds on the rate of entropy production by
  antidiffusive fluxes and Runge--Kutta (RK) time discretizations.
  Two versions of the fully discrete entropy fix are developed
  for this purpose. The first one incorporates temporal entropy
  production into the flux constraints, which makes them
  more restrictive and dependent on the time step. The second
  algorithm interprets the final stage of a high-order AFC-RK
  method as a constrained antidiffusive correction of an implicit
  low-order scheme (algebraic Lax--Friedrichs in space + backward
  Euler in time). In this case, iterative flux correction is
  required, but the inequality constraints are less restrictive
  and limiting can be performed using algorithms developed
  for the semi-discrete problem. To motivate the use of limiter-based
  entropy fixes, we prove a finite element version of the
  Lax--Wendroff theorem and perform numerical studies for
  standard test problems. In our numerical experiments,
  entropy-dissipative schemes converge to correct
  weak solutions of scalar conservation laws, of
  the Euler equations,  and of the shallow water equations.
  
\end{abstract}

\begin{keyword}
  hyperbolic conservation laws \sep property-preserving schemes\sep
  continuous finite elements 
 \sep algebraic flux correction\sep
  convex limiting \sep entropy stabilization

\end{keyword}
\end{frontmatter}

\section{Introduction}
It is well known that nonlinear hyperbolic problems may have multiple weak solutions but the one corresponding to the vanishing viscosity limit is unique and satisfies a weak form of an entropy inequality. Moreover, the conserved variables or derived quantities are often known to be bounded in a certain manner. It is therefore essential to use numerical methods that are bound preserving and entropy dissipative. Failure to do so may cause occurrence of nonphysical states or convergence to wrong weak solutions \cite{banks2009,chandrashekar2013,einfeldt1988,kurganov2007,toro2009}.

Many modern high-resolution schemes use limiters to ensure preservation of local bounds or at least positivity preservation
for scalar quantities of interest. In the context of finite element approximations, such schemes can be constructed using geometric slope limiting \cite{zhang2011} or the framework of algebraic flux correction (AFC) \cite{kuzmin2012a} and its extensions to hyperbolic systems \cite{guermond2018,hajduk2021,kuzmin2012b,kuzmin2020}. The additional requirement of entropy stability implies that a \mbox{(semi-)discrete} entropy inequality should hold for at least one entropy pair. Entropy-dissipative space discretizations of second and higher order can be designed as in \cite{chen2017,ray2016,kuzmin2020f} using Tadmor's criterion \cite{tadmor2003} of entropy stability for semi-discrete problems. Recent years have witnessed significant advances in the development of high-order entropy-stable schemes that exploit the summation-by-parts (SBP) property of discrete operators collocated at quadrature points \cite{carpenter2014,chen2020,fisher2013,gassner2013,pazner2019}. A way to penalize semi-discrete entropy production by solution gradients inside mesh cells was proposed by Abgrall \cite{abgrall2018} and generalized in \cite{abgrall2019,kuzmin2020g}. To ensure global entropy stability of a fully discrete scheme, relaxation Runge--Kutta methods were developed by Ketcheson et al. \cite{ketcheson2019,ranocha2020}. The fully discrete TECNO-SSP schemes proposed in \cite{zakerzadeh2015} are globally entropy stable under mild time step restrictions. Methods that ensure the validity of local fully discrete entropy inequalities are still rare. Their development was recently advanced by the work of Kivva \cite{kivva-arxiv} and Berthon et al.\ \cite{berthon2021,berthon2020}. The algorithms developed in these publications are closely related to our approach. However, they are currently restricted to explicit finite volume schemes on uniform meshes in 1D.

As of this writing, very few bound-preserving schemes for nonlinear
hyperbolic problems are provably entropy stable and vice versa. The
AFC scheme proposed in \cite{kuzmin2020f} constrains a continuous
finite element discretization of a scalar conservation law using
a bound-preserving flux limiter and a semi-discrete entropy fix
based on Tadmor's condition. In the present paper, we extend the
underlying methodology to arbitrary fluxes, hyperbolic systems,
and fully discrete schemes. Instead of constructing entropy-conservative
fluxes and adding second-order entropy viscosity, we limit the
antidiffusive flux that transforms a property-preserving low-order
method into a high-order baseline scheme. We show that the bounds of the 
inequality constraints for limiter-based entropy fixes can be defined to
produce an entropy-conservative or entropy-dissipative scheme. We
derive sufficient conditions for the validity of local entropy
inequalities and constrain the antidiffusive fluxes accordingly.
In addition to a generalized entropy fix for the spatial
semi-discretization, we present two algorithms that ensure fully
discrete entropy stability. The first one limits entropy
production by forward
Euler stages of a strong stability preserving Runge--Kutta (RK)
method, while the second one constrains the final
stage of a general RK method in an iterative manner. We
study the effectiveness of each limiter-based entropy fix numerically and prove
a Lax--Wendroff-type theorem about convergence of finite element
approximations to entropy solutions.

In Section~\ref{sec:fem} of this paper, we discretize a generic
hyperbolic problem (a scalar conservation law or a system) using
the continuous Galerkin method and (multi-)linear finite
elements. Next, we review the basic principles of algebraic flux
correction for space discretizations of this kind. Existing
extensions to high-order finite elements and discontinuous
Galerkin (DG) methods are cited as well.
The new limiter-based entropy fixes are presented 
in Section~\ref{sec:efix}. The importance of
entropy stability for convergence to vanishing viscosity solutions is
illustrated by the results of theoretical studies in Section~\ref{sec:lw}
and numerical experiments in Section~\ref{sec:results}. We close
this paper with concluding remarks in Section~\ref{sec:conclusion}.

\section{Continuous FEM for hyperbolic problems}
\label{sec:fem}

Let $u(\vec{x},t)\in\R^{m\times 1}$ denote the local density of $m\in\mathbb{N}$
conserved quantities at the space location $\vec{x}\in\R^d, d\in\{1,2,3\}$
and time $t\ge 0$. Imposing periodic boundary conditions on the boundaries
of a spatial domain $\Omega\subset\R^d$, we consider the initial 
value problem
\begin{subequations}\label{ibvp}
\begin{align}
\pd{u}{t}+\nabla\cdot\vec{f}(u)&=0\quad \mbox{in}\ \Omega\times\R^+,
\label{ibvp-pde}\\
\qquad u(\cdot,0)&=u_0\quad\mbox{in}\ \Omega,\label{ibvp-ic}
\end{align}
\end{subequations}
where $\vec{f}(u)=(f_{ki})
\in\R^{m\times d}$
is an array of inviscid fluxes and $\nabla\cdot\vec{f}=
\left(\sum_{i=1}^d\pd{f_{ki}}{x_i}\right)\in\R^{m\times 1}$.
The flux and Jacobian of a projection onto a 
vector $\mathbf{n}=(n_i)\in\R^{d\times 1}$
are defined by 
$$\vec f(u)\cdot\vec n=\left(
\sum\limits_{i=1}^df_{ki}n_i\right)\in\mathcal\R^{m\times 1},\qquad
\vec f'(u)\cdot\vec n=\left(\sum\limits_{i=1}^d
\pd{f_{ki}}{u_l}n_i\right)\in\R^{m\times m}.$$
In the multidimensional case ($d>1$), the
 array \smash{\mbox{$\vec f'(u)=\left(\pd{f_{ki}}{u_l}\right)
\in\R^{m\times m\times d}$}} is composed from Jacobian matrices
$\vec f'(u)\cdot\vec e_1,\ldots,\vec f'(u)\cdot\vec e_d$,
where $\vec e_i$ is the $i$-th unit vector in $\R^d$.
We assume that \eqref{ibvp-pde}
is hyperbolic, i.e., that any directional Jacobian 
$\vec{f}'(u)\cdot\mathbf{n}$
is diagonalizable
with real eigenvalues $\lambda_1,\ldots,\lambda_m$
representing $m$ finite speeds of wave propagation.

Because of hyperbolicity, there exists a convex
entropy $\eta:\R^m\to\R$ and associated entropy flux
$\vec{q}:\R^m\to\R^{1\times d}$ such that 
$\vec q'(u)=\eta'(u)^T\vec f'(u)$ for $\eta'(u)
=\left(\pd{\eta}{u_k}\right)\in\R^{m\times 1}$. If problem \eqref{ibvp}
has a smooth classical solution $u$, the entropy conservation law
\beq\label{ent-cons}
 \pd{\eta(u)}{t}+\nabla\cdot\vec{q}(u)=0 \qquad\mbox{in}\ \Omega\times\R^+
 \eeq
 can be derived from \eqref{ibvp-pde} using multiplication by the
 vector $v(u)=\eta'(u)$ of entropy
 variables, the chain rule, and the definition of an
 entropy pair $\{\eta(u),\vec{q}(u)\}$. In general, the vanishing
 viscosity solution of \eqref{ibvp-pde} satisfies a weak form
 of the entropy inequality
\beq
 \pd{\eta(u)}{t}+\nabla\cdot\vec{q}(u)\le 0 \qquad\mbox{in}\ \Omega\times\R^+
\label{ent-ineq}
\eeq
for any entropy pair. Hence, entropy is conserved in
smooth regions and dissipated at shocks. For the derivation
of \eqref{ent-ineq} and further details, we refer the
reader to Tadmor \cite{tadmor2003}.

When it comes to solving
\eqref{ibvp-pde} numerically, it is essential to guarantee
that a (semi-) discrete version of \eqref{ent-ineq} holds for at
least one entropy pair. Moreover, a  good numerical
method should ensure preservation of invariant domains, i.e.,
produce approximations belonging to a convex set $\mathcal G$ if the
exact solution of \eqref{ibvp-pde} 
is known to stay in this set \cite{guermond2016}.

\subsection{Consistent Galerkin discretization}

We discretize \eqref{ibvp-pde} in space using the continuous
Galerkin method on a conforming mesh $\mathcal T_h=\{
K^{1},\ldots,K^{E_h}\}$ of linear
($\mathbb{P}_1$) or multilinear ($\mathbb{Q}_1$) finite elements. 
The numerical solution $u_h=\sum_{j=1}^{N_h}u_j\varphi_j$ is
expressed in terms of Lagrange basis functions 
$\varphi_1,\ldots,\varphi_{N_h}$ such that
$\sum_{i=1}^{N_h}\varphi_i\equiv 1$. The basis
function $\varphi_i$ is associated with a vertex
$\vec{x}_i$ of $\mathcal T_h$. It has the property
that  $\varphi_i(\vec{x}_j)=\delta_{ij}$.
The indices of
elements containing $\vec{x}_i$ are stored in
the integer set $\mathcal E_i$. The global indices
of nodes belonging to $K^e\in\mathcal T_h$ are stored
in the integer set $\mathcal N^e$. The computational
stencil of node $i$ is defined by the index set
$\mathcal N_i=\bigcup_{e\in\mathcal E_i}\mathcal N^e$.

The semi-discrete weak form of \eqref{ibvp-pde}
with periodic boundary conditions is given by
\beq\label{semi-gal1}
\sum_{e=1}^{E_h}\int_{K^e}v_h\left[
\pd{u_h}{t}+\nabla\cdot\vec f(u_h)\right]\dx=0\qquad
\forall v_h\in V_h,
\eeq
where $V_h$ is the finite-dimensional space spanned by the
basis functions $\varphi_1,\ldots,\varphi_{N_h}$. For any
entropy pair $\{\eta(u),\mathbf{q}(u)\}$ and $v(u)=\eta'(u)$,
the evolution of $\eta(u_h)$ is governed by \cite{kuzmin2020g}
  \beq\label{cons-eta-gal}
 \sum_{e=1}^{E_h}\int_{K^e}\left[\pd{\eta(u_h)}{t}+
\nabla\cdot\mathbf{q}(u_h)\right]\dx=\sum_{e=1}^{E_h}
\int_{K^e}v(u_h)\left[\pd{u_h}{t}+\nabla\cdot\mathbf{f}(u_h)\right]\dx.
\eeq
In the scalar case ($m=1$), the right-hand side of
\eqref{cons-eta-gal} vanishes due to \eqref{semi-gal1} for
$v(u_h)=u_h$ corresponding to $\eta(u)=\frac{u^2}2$.
It follows that the continuous Galerkin discretization
of \eqref{ibvp-pde}  satisfies
an integral form of \eqref{ent-cons}, i.e.,
is entropy conservative 
in this particular case \cite{kuzmin2020g}. This remarkable property was first discovered
by Tadmor \cite{tadmor1986,tadmor2003} in the context of finite volume
schemes written as lumped-mass
$\mathbb{P}_1$ finite element approximations.

Invoking
the definition of $u_h$ and using the test function
$v_h=\varphi_i$ in \eqref{semi-gal1}, we obtain
\beq\label{semi-gal2}
\sum_{j\in\mathcal N_i}m_{ij}\td{u_j}{t}+
\sum_{e\in\mathcal E_i}\int_{K^e}\varphi_i
\nabla\cdot\vec f(u_h)\dx=0,
\eeq
where
$$
m_{ij}=\sum_{e\in\mathcal E_i\cap\mathcal E_j}\int_{K^e}\varphi_i
\varphi_j\dx
$$
is an entry of the consistent mass matrix $M_C=(m_{ij})_{i,j=1}^{N_h}$. Since
$u_h$ is differentiable on $K^e$,
equation \eqref{semi-gal1} can be written in the equivalent
quasi-linear form
\beq\label{semi-gal3}
\sum_{j\in\mathcal N_i}\left[m_{ij}\td{u_j}{t}+a_{ij}(u_h)u_j\right]=0,
\eeq
where
$$
a_{ij}(u_h)=\sum_{e\in\mathcal E_i\cap\mathcal E_j}\int_{K^e}\varphi_i
\vec f'(u_h)\cdot\nabla \varphi_j\dx
$$
is an entry of the solution-dependent discrete Jacobian operator $A(u_h)=(a_{ij}(u_h))_{i,j=1}^{N_h}$.

\subsection{Quadrature-based approximations}

The properties of a finite element discretization
may change if the integrals are calculated
numerically. For example, the use of inexact nodal
quadrature transforms the matrix $M_C$
into its lumped counterpart
$
M_L=(\delta_{ij}m_i)_{i,j=1}^{N_h}
$ with diagonal entries
$$
m_i=
\sum_{e\in\mathcal E_i}\int_{K^e}\varphi_i\dx
=
\sum_{e\in\mathcal E_i}\int_{K^e}\varphi_i
\underbrace{\left(\sum_{j\in\mathcal N_i}\varphi_j\right)}_{\equiv 1}\dx=
\sum_{j\in\mathcal N_i}m_{ij}.
$$

The second-order accurate group finite element approximation
\cite{barrenechea2017b,fletcher1983,kuzmin2010a}
\beq\label{groupflux}\vec f(u_h)\approx
\mathbf{f}_h(u_h)=\sum_{j=1}^{N_h}\mathbf{f}_j\,\varphi_j,\qquad
\mathbf{f}_j=\mathbf{f}(u_j)
\eeq
yields a useful inexact quadrature rule for
calculating \smash{\mbox{$\int_{K^e}\varphi_i\nabla\cdot\vec f(u_h)\dx$}}.
Replacing $\vec f(u_h)$ by
$\mathbf{f}_h(u_h)$ in the lumped-mass version of
\eqref{semi-gal2}, we obtain the semi-discrete scheme
\beq\label{semi-group}
m_i\td{u_i}{t}+
\sum_{j\in\mathcal N_i}\mathbf{f}_j\cdot\bc_{ij}=0,
\eeq
where
$$
\bc_{ij}=\sum_{e\in\mathcal E_i\cap\mathcal E_j}\int_{K^e}\varphi_i
\nabla \varphi_j\dx
$$
denotes a vector-valued entry of the discrete gradient operator
$\mathbf{C}=(\bc_{ij})_{i,j=1}^{N_h}$. Note that
\smash{\mbox{$\sum_{j\in\mathcal N_i}\bc_{ij}=\sum_{j=1}^{N_h}\bc_{ij}=
\mathbf{0}$}} since \smash{\mbox{$\sum_{j=1}^{N_h}\nabla\varphi_j
\equiv\nabla 1=\mathbf{0}$}} for Lagrange basis functions $\varphi_j$.
Moreover, $\bc_{ji}=-\bc_{ij}$ unless $i$ or $j$ is
a node on a non-periodic boundary. It follows that
\beq\label{props-cij}
\sum_{j\in\mathcal N_i}\mathbf{f}_j\cdot\bc_{ij}
=\sum_{j\in\mathcal N_i\backslash\{i\}}(\bF_j-\bF_i)\cdot\bc_{ij}
=\sum_{j\in\mathcal N_i\backslash\{i\}}(\bF_j+\bF_i)\cdot\bc_{ij}
=\sum_{j\in\mathcal N_i}(\bF_j+\bF_i)\cdot\bc_{ij}.
\eeq
As shown in \cite{selmin1993,selmin1996}, the quadrature-based
approximation \eqref{semi-group} is equivalent to a
vertex-centered finite volume method with a centered
numerical flux. Many edge-based finite element schemes
for compressible flow problems exploit this relationship
to finite volumes when it comes to the design of artificial
viscosities and flux limiting \cite{kuzmin2012a,lyra2002,lohner2008}.

\begin{rmk}
Similarly to DG schemes with inexact quadrature \cite{chen2017}, the
use of the group finite element formulation \eqref{groupflux} may
cause a lack of entropy stability even in the case $m=1,$ $\eta(u)
=\frac{u^2}2$, in which
the consistent Galerkin approximation is entropy stable.
  \end{rmk}

\subsection{Algebraic flux correction}

The concept of algebraic flux correction (AFC)
provides a general framework for converting \eqref{semi-group}
into a  property-preserving space discretization of the form
\cite{kuzmin2012a,kuzmin2012b,lohmann2019}
\beq\label{semi-afc}
m_i\td{u_i}{t}=
\sum_{j\in\mathcal N_i\backslash\{i\}}
    [d_{ij}(u_j-u_i)-(\mathbf{f}_j-\mathbf{f}_i)
      \cdot\bc_{ij}+f_{ij}^*].
\eeq
The coefficients of the artificial viscosity
operator $D=(d_{ij})_{i,j=1}^{N_h}$ have the property that
$$
\sum_{j\in\mathcal N_i}d_{ij}=0,\qquad
d_{ij}=d_{ji},\qquad  i,j=1,\ldots,N_h
$$
and are chosen in such a way that all relevant inequality
constraints (maximum principles, entropy stability
conditions etc.) are satisfied in the case $f_{ij}^*=0$.
We define such artificial viscosity coefficients $d_{ij}$
in Section \ref{sec:alf}. In general, 
$f_{ij}^*=-f_{ji}^*$ is supposed to approximate a target
flux $f_{ij}=-f_{ji}$ in a property-preserving manner. We
discuss the corresponding limiting procedures in Sections
\ref{sec:mcl}, \ref{sec:timelim}, and \ref{sec:semi}--\ref{sec:imp}.
The discretization to which \eqref{semi-afc} reduces in
the case $f_{ij}^*=f_{ij}\
\forall j\in\mathcal N_i\backslash\{i\}$ is typically a
stabilized version of \eqref{semi-gal2}. It may use the
group finite element approximation~\eqref{groupflux} and/or
high-order stabilization provided that the order of
accuracy is preserved. Definitions of (stabilized)
target fluxes $f_{ij}$ for $\mathbb{P}_1/\mathbb{Q}_1$ and
higher-order finite element discretizations of \eqref{ibvp-pde}
can be found in \cite{kuzmin2020,kuzmin2020f,kuzmin2020g}. The
fluxes
that we constrain in the numerical experiments of this work are
defined in Section \ref{sec:results}.

\section{Invariant domain preserving schemes}
\label{sec:idp}

\subsection{Algebraic Lax--Friedrichs method}
\label{sec:alf}

As explained above, the definition of the scalar artificial
viscosity coefficient $d_{ij}$ for the flux-corrected
scheme \eqref{semi-afc} should ensure that its
low-order counterpart
\beq\label{semi-low}
m_i\td{u_i}{t}=
\sum_{j\in\mathcal N_i\backslash\{i\}}
    [d_{ij}(u_j-u_i)-(\mathbf{f}_j-\mathbf{f}_i)
      \cdot\bc_{ij}]
    \eeq
 is property preserving.
 An algebraic version of the local Lax--Friedrichs method uses
 \cite{guermond2016,kuzmin2020}
    \beq
d_{ij} =\begin{cases}
\max\{\lambda_{ij}|\mathbf{c}_{ij}|,
\lambda_{ji}|\mathbf{c}_{ji}|\} & \mbox{if}\ j\ne i,\\
  -\sum_{k\in \mathcal N_i\backslash\{i\}}  d_{ik}
  & \mbox{if}\ j=i,
\end{cases}
  \label{gms}
  \eeq
  where $\lambda_{ij}$ is an upper bound for the
  spectral radius of
  $A_{ij}(\hat u)=\mathbf{f}'(\hat u)\cdot\frac{\bc_{ij}}{|\bc_{ij}|}$
  evaluated at an arbitrary convex combination $\hat u$ of the
  states $u_i$ and $u_j$.
  As shown by Guermond and Popov \cite{guermond2016},
  this definition of the maximum wave speed
  $\lambda_{ij}$ guarantees preservation
  of invariant domains and entropy stability. The first
  property is a generalized maximum principle. If there
  is a convex invariant set $\mathcal G$ such that
  $u_h(\bx,0)\in\mathcal G$ a.e.\ in $\Omega$, an invariant
  domain preserving (IDP) scheme keeps $u_h(\bx,t)$ in $\mathcal G$
  for all $t>0$. The entropy stability property is needed
  to avoid convergence to wrong weak solutions. Entropy
  stable schemes produce approximations 
  that satisfy a (semi-)discrete entropy
  inequality consistent with \eqref{ent-ineq}.
  
\begin{rmk}
  In earlier versions  \cite{kuzmin2012b,kuzmin2010a,lyra2002} of
  the algebraic Lax--Friedrichs (ALF) scheme
\eqref{semi-low},\eqref{gms} for nonlinear
  problems, $\lambda_{ij}$ was taken to be 
  the spectral radius  of the 
  Jacobian $A_{ij}(u_i)$ or of the  Roe matrix
  $\hat A_{ij}\in\R^{m\times m}$ such that $\hat A_{ij}(u_j-u_i)
  =(\mathbf{f}_j-\mathbf{f}_i)\cdot\bc_{ij}$.
  Although these definitions work well in practice,
  the resulting schemes are not
  provably property-preserving for systems. For
  linear advection of a scalar conserved quantity
  with constant velocity, all versions of the
  ALF method reduce to the discrete upwinding procedure
  described in \cite{kuzmin2012a}.
\end{rmk}
  
\subsection{Bound-preserving convex limiting}
\label{sec:mcl}

To formulate a sufficient condition for an AFC scheme
of the form \eqref{semi-afc} to inherit
the IDP property of the ALF method \eqref{semi-low},
we consider the equivalent representation
\beq\label{semi-mcl}
m_{i}\td{u_i}{t}=\sum_{j\in\mathcal N_i\backslash\{i\}}
2d_{ij}(\bar u_{ij}^*-u_i),\qquad \bar u_{ij}^*=
\bar u_{ij}+\frac{f_{ij}^*}{2d_{ij}}
\eeq
of \eqref{semi-afc} in terms of the property-preserving
intermediate states
\beq\label{barstate}
\bar u_{ij}=\frac{u_j+u_i}{2}-\frac{(\mathbf{f}_j-\mathbf{f}_i)\cdot\mathbf{c}_{ij}}{2d_{ij}}.
\eeq
Let $\mathcal G$ be a convex invariant set of the initial
value problem \eqref{ibvp}.
As noticed by Guermond and Popov \cite{guermond2016},
 the ALF \emph{bar state} $\bar u_{ij}$ represents
 an averaged exact solution of a one-dimensional
 Riemann problem with
the initial states $u_i$ and $u_j$. It follows that
$\bar u_{ij}\in\mathcal G$ whenever $u_i,u_j\in\mathcal G$.
Hence, a space discretization of the form \eqref{semi-mcl}
is IDP if $f_{ij}^*$ is limited in such a way that
$\bar u_{ij}^*\in\mathcal G$ for $\bar u_{ij}\in\mathcal G$. 
The monolithic convex limiting (MCL) algorithm proposed
in \cite{kuzmin2020} is designed to ensure this property
and validity of local maximum principles for scalar
quantities of interest. Extensions of the MCL methodology
to high-order Bernstein finite elements and discontinuous Galerkin methods can be found in \cite{hajduk2021,kuzmin2020g}.

\subsection{Bound-preserving time integration}
\label{sec:timelim}

If the system of semi-discrete equations \eqref{semi-mcl}
is integrated in time using an explicit strong stability
preserving (SSP) Runge--Kutta method, each Shu--Osher stage
is of the form
\beq\label{ssp-mcl}
u_i^{\mathrm{SSP}}
=u_i+\frac{\Delta t}{m_i}\sum_{j\in\mathcal N_i\backslash\{i\}}
2d_{ij}(\bar u_{ij}^*-u_i)
=(1-c_i)u_i+c_i\bar u_i^*,
\eeq
where
$$
c_i=\frac{\Delta t}{m_i}\sum_{j\in\mathcal N_i\backslash\{i\}}
2d_{ij},\qquad \bar u_i^*=\frac{1}{\sum_{j\in\mathcal N_i\backslash\{i\}}
2d_{ij}}\sum_{j\in\mathcal N_i\backslash\{i\}}
2d_{ij}\bar u_{ij}^*.
$$
Note that $\bar u_i^*$ is a convex combination of $\bar u_{ij}^*,
\ j\in\mathcal N_i\backslash\{i\}$.
If the time step $\Delta t$ satisfies the CFL-like condition
$c_i\le 1$, then $u_i^{\mathrm{SSP}}$ is a convex combination
of $u_i$ and $\bar u_i^*$. Thus
$$
u_j\in\mathcal G\ \forall j\in\mathcal N_i\
\implies\ \bar u_{ij}^*\in\mathcal G\ \forall j\in\mathcal N_i\backslash\{i\}\ \implies\ \bar u_i^*\in\mathcal G
\ \implies\ u_i^{\mathrm{SSP}}\in\mathcal G.
$$
The final stage of a general high-order Runge--Kutta method
for \eqref{semi-mcl} can also be written in the form
\eqref{ssp-mcl} and constrained to produce
$\bar u_i^*\in\mathcal G$; see \cite{kuzmin-arxiv} for details.

\begin{rmk}
The IDP property of the fully discrete scheme
can also be enforced using convex limiting of predictor-corrector
type \cite{guermond2018,pazner2021}. However, such flux-corrected
transport (FCT) algorithms are not well suited for entropy
fixes that we propose in Section \ref{sec:efix} because
the underlying theory requires monolithic flux correction
at the semi-discrete level.
\end{rmk}

\section{Entropy stabilization via limiting}
\label{sec:efix}

Since the low-order ALF method \eqref{semi-low} is
property-preserving, entropy stability of the AFC scheme
\eqref{semi-afc} can always be enforced by further reducing the
magnitude of the fluxes $f_{ij}^*$ if necessary. In this
section, we extend the limiter-based entropy correction
techniques proposed in \cite{kuzmin2020f,kuzmin2020g} to
arbitrary target fluxes and fully discrete nonlinear schemes.

\subsection{Semi-discrete entropy correction}
\label{sec:semi}

An entropy-corrected
version of the semi-discrete AFC scheme
\eqref{semi-afc} is defined by
\begin{align}\label{semi-afc-fix}
m_i\td{u_i}{t}=
\sum_{j\in\mathcal N_i\backslash\{i\}}
    [d_{ij}(u_j-u_i)-(\mathbf{f}_j-\mathbf{f}_i)
      \cdot\bc_{ij}+\alpha_{ij}f_{ij}^*].
\end{align}
To make it entropy conservative/dissipative w.r.t. an
entropy pair $\{\eta(u),\mathbf{q}(u)\}$,
we apply correction factors
$\alpha_{ij}\in[0,1]$ satisfying $\alpha_{ij}=\alpha_{ji}$ and
Tadmor's condition \cite{tadmor2003}
\beq\label{estab}
\frac{(v_i-v_j)^T}2[
d_{ij}(u_j-u_i)-
       (\mathbf{f}_j+\mathbf{f}_i) \cdot\mathbf{c}_{ij}
+\alpha_{ij}f_{ij}^*]\le (\boldsymbol{\psi}_j
  -\boldsymbol{\psi}_i) \cdot\mathbf{c}_{ij},
\eeq
where $v_i=\eta'(u_i)$ and $\boldsymbol{\psi}_j
=v_j^T\mathbf{f}_j-\mathbf{q}_j$ is the entropy potential
corresponding to $\mathbf{q}_j=\mathbf{q}(u_j)$.

Following the proofs in \cite{ray2016} and
\cite{kuzmin2020f}, we multiply
\eqref{semi-afc-fix} by $v_i^T$ and use the zero sum property
$\sum_{j\in\mathcal N_i}\bc_{ij}=\mathbf{0}$ of the
discrete gradient operator to show that \eqref{estab} implies
\begin{align}
  m_iv_i^T\td{u_i}{t}&=
\sum_{j\in\mathcal N_i\backslash\{i\}}
   v_i^T[d_{ij}(u_j-u_i)-(\mathbf{f}_j-\mathbf{f}_i)
     \cdot\bc_{ij}+\alpha_{ij}f_{ij}^*]\nonumber\\
  &=
\sum_{j\in\mathcal N_i\backslash\{i\}}
   v_i^T[d_{ij}(u_j-u_i)-(\mathbf{f}_j+\mathbf{f}_i)
     \cdot\bc_{ij}+\alpha_{ij}f_{ij}^*]-2
   v_i^T\mathbf{f}_i\cdot\mathbf{c}_{ii}
   \nonumber\\
   &=\sum_{j\in\mathcal N_i\backslash\{i\}}
   \frac{(v_i+v_j)^T}2[d_{ij}(u_j-u_i)-(\mathbf{f}_j+\mathbf{f}_i)
     \cdot\bc_{ij}+\alpha_{ij}f_{ij}^*]\nonumber\\
   &+\sum_{j\in\mathcal N_i\backslash\{i\}}
   \frac{(v_i-v_j)^T}2[d_{ij}(u_j-u_i)-(\mathbf{f}_j+\mathbf{f}_i)
     \cdot\bc_{ij}+\alpha_{ij}f_{ij}^*]-2
   v_i^T\mathbf{f}_i\cdot\mathbf{c}_{ii}\nonumber\\
&\le\sum_{j\in\mathcal N_i\backslash\{i\}}
   \frac{(v_i+v_j)^T}2[d_{ij}(u_j-u_i)-(\mathbf{f}_j+\mathbf{f}_i)
     \cdot\bc_{ij}+\alpha_{ij}f_{ij}^*]\nonumber\\
&+\sum_{j\in\mathcal N_i\backslash\{i\}}(\boldsymbol{\psi}_j
     -\boldsymbol{\psi}_i) \cdot\mathbf{c}_{ij}
     -2v_i^T\mathbf{f}_i\cdot\mathbf{c}_{ii}\nonumber\\
     &=\sum_{j\in\mathcal N_i\backslash\{i\}}
   \frac{(v_i+v_j)^T}2[d_{ij}(u_j-u_i)-(\mathbf{f}_j+\mathbf{f}_i)
     \cdot\bc_{ij}+\alpha_{ij}f_{ij}^*]\nonumber\\
   &+\sum_{j\in\mathcal N_i\backslash\{i\}}
   [v_j^T\mathbf{f}_j+v_i^T\mathbf{f}_i-
   (\mathbf{q}_j-\mathbf{q}_i)]\cdot\mathbf{c}_{ij}\nonumber\\
   &=\sum_{j\in\mathcal N_i\backslash\{i\}}[G_{ij}-
   (\mathbf{q}_j-\mathbf{q}_i)\cdot\bc_{ij}].\label{semi-ineq}
\end{align}  
The fluxes $G_{ij}$ that appear on the right-hand side of
the last equation are given by \cite{kuzmin2020f}
$$
 G_{ij}=\frac{(v_i+v_j)^T}2[d_{ij}(u_j-u_i)+\alpha_{ij}f_{ij}^*]
-\frac{(v_i-v_j)^T}2(\mathbf{f}_j-\mathbf{f}_i)\cdot\bc_{ij}.
$$
Since $v_i^T\td{u_i}{t}=\td{\eta(u_i)}{t}$ by
definition of the entropy variable
$v_i$ and the chain rule, the semi-discrete entropy
inequality \eqref{semi-ineq} represents a consistent
space discretization of \eqref{ent-ineq}.

\begin{rmk}
 In the case of periodic boundary conditions, we have
$\mathbf{c}_{ji}=-\mathbf{c}_{ij}$ and therefore
$G_{ji}=-G_{ij}$ for all $j\in\mathcal N_i\backslash\{i\}$.
Moreover, periodicity implies that
\beq
\sum_{i=1}^{N_h}\sum_{j\in\mathcal N_i\backslash\{i\}}[G_{ij}-
   (\mathbf{q}_j-\mathbf{q}_i)\cdot\bc_{ij}]=
\sum_{i=1}^{N_h}
\sum_{j\in\mathcal N_i\backslash\{i\}}[G_{ij}-
  (\mathbf{q}_j+\mathbf{q}_i)\cdot\bc_{ij}]=0.
\eeq
Although we do assume periodicity for ease of presentation,
the below fix is derived for the general case and
can be used for element-level entropy corrections (as in
\cite{kuzmin2020f,kuzmin2020g}).
\end{rmk}

Let us now discuss the practical calculation of  $\alpha_{ij}
\in[0,1]$.
Condition \eqref{estab} holds if the rate of
entropy production $\frac12(v_i-v_j)^T\alpha_{ij}f_{ij}^*$
by the flux $\alpha_{ij}f_{ij}^*$
does not exceed
\beq\label{Qij-EC}
Q_{ij}\le Q_{ij}^{\mathrm{EC}}:=(\boldsymbol{\psi}_j
-\boldsymbol{\psi}_i) \cdot\mathbf{c}_{ij}
-\frac{(v_i-v_j)^T}2[d_{ij}(u_j-u_i)-
       (\mathbf{f}_j+\mathbf{f}_i) \cdot\mathbf{c}_{ij}].
\eeq

As shown by Chen and Shu \cite{chen2017} in the context of
DG methods, the stability condition \eqref{estab} is
satisfied for $\alpha_{ij}=0$. Hence, the entropy
conservative upper bound $Q_{ij}^{\mathrm{EC}}$ is nonnegative.

The limiter-based entropy correction procedure
proposed in \cite{kuzmin2020f}
multiplies $f_{ij}^*$ by 
\beq\label{alpha}
\alpha_{ij}=
\begin{cases}
\displaystyle
  \min\left\{1,
    \frac{\min\{Q_{ij},\frac12(v_i-v_j)^Tf_{ij}^*,Q_{ji}\}}{\frac12(v_i-v_j)^Tf_{ij}^*}\right\}
  & \mbox{if}\ (v_i-v_j)^Tf_{ij}^*>0,\\ 
1 &\mbox{otherwise}
\end{cases}
\eeq
with $Q_{ij}=Q_{ij}^{\mathrm{EC}}$ and adds second-order entropy
viscosity to the target flux $f_{ij}$ before IDP limiting. In
the present work, we leave $f_{ij}$ unchanged but impose 
an entropy-dissipative bound $Q_{ij}\le Q_{ij}^{\mathrm{EC}}$
on the rate of entropy production by $\alpha_{ij}f_{ij}^*$.
Adapting the 
entropy viscosity formula that was used to stabilize
the target fluxes in \cite{kuzmin2020f}, we define
\beq\label{Qij-ED}
Q_{ij}^{\mathrm{ED}}=\max\left\{0,
Q_{ij}^{\mathrm{EC}}+\min\left\{0,\frac{(v_i-v_j)^T}2\left[
  \mathbf{f}_j+\mathbf{f}_i-2\mathbf{f}\left(
  \frac{u_j+u_i}2\right)\right]\cdot\bc_{ij}\right\}
\right\}
\eeq
and use $Q_{ij}=Q_{ij}^{\mathrm{ED}}$ in \eqref{alpha} to generate
sufficient amounts of entropy viscosity 
at shocks.
\begin{rmk}
  The limiter-based approach makes it possible to
  satisfy \eqref{estab}  for multiple entropies
  by using the minimum of the corresponding correction factors
  $\alpha_{ij}$.
 \end{rmk}
\begin{rmk}
  In practice, we calculate the entropy correction factors $\alpha_{ij}$
  as follows:
  $$
  \alpha_{ij}=
\begin{cases}\displaystyle
  \frac{2\min\{Q_{ij},Q_{ji}\}+\delta|f_{ij}^*|}{
    (v_i-v_j)^Tf_{ij}^*+\delta|f_{ij}^*|}
  & \mbox{if}\ (v_i-v_j)^Tf_{ij}^*>2\min\{Q_{ij},Q_{ji}\},\\ 
1 &\mbox{otherwise},
\end{cases}$$
where $\delta$ is a small positive number (we use $\delta=10^{-2}$ in
the numerical experiments of Section~\ref{sec:results}). This
definition  produces
$\alpha_{ij}=1$ in the limits $|v_i-v_j|\to 0$
and $|f_{ij}^*|\to 0$. Importantly,
it ensures continuity of $\alpha_{ij}f_{ij}^*$
which is needed for well-posedness of
 fully discrete nonlinear problems
 and convergence
 of iterative solvers for implicit schemes.
\end{rmk}

\subsection{Fully discrete explicit correction}
\label{sec:exp}

The entropy fixes presented so far are designed to ensure that
the semi-discrete AFC scheme is entropy stable. However,
discretization in time may increase the rate of entropy
production, leading to a fully discrete scheme that lacks
entropy stability \cite{merriam1989,tadmor2003}. While the fully implicit backward Euler
discretization of \eqref{semi-afc-fix} is property preserving
for any time step, it is only first-order
accurate in time and requires iterative solution of nonlinear systems.
As shown by Guermond and Popov \cite{guermond2016},
forward Euler stages \eqref{ssp-mcl} of
an explicit SSP Runge--Kutta method are  entropy
stable in the case $\bar u_{ij}^*:=\bar u_{ij}$, i.e.,
if all fluxes $f_{ij}^*$ are set to zero. Nonvanishing
 fluxes
$\alpha_{ij}f_{ij}^*$ may produce too much entropy even if
they satisfy Tadmor's stability condition \eqref{estab}.
In this section, we limit them in a way that ensures fully
discrete entropy stability under suitable assumptions.
In the next section, we present an iterative fix based
on \eqref{estab}.
For further discussion and analysis of fully discrete
entropy stability, we refer the interested reader to
LeFloch et al. \cite{lefloch2002}, Lozano \cite{lozano2018,lozano2019},
and Merriam \cite{merriam1989}.

Each forward Euler stage of an explicit SSP-RK method for
\eqref{semi-afc-fix} can be written as
\beq
u_i^{\mathrm{SSP}}=u_i+\Delta t\dot u_i^*,
\eeq
$$
\dot u_i^*=\frac{1}{m_i}\sum_{j\in\mathcal N_i\backslash\{i\}}
     [d_{ij}(u_j-u_i)-
       (\mathbf{f}_j-\mathbf{f}_i) \cdot\mathbf{c}_{ij}
    +\alpha_{ij}f_{ij}^*].
$$
Using Taylor's theorem, we find that (cf. \cite{merriam1989}, Section 4.2)
 \beq\label{ent-taylor}
    \frac{\eta(u_i^{\mathrm{SSP}})-\eta(u_i)}{\Delta t}
    =v_i^T\dot u_i^*
       + \frac{\Delta t}{2}
      (\dot u_i^*)^T\eta''(\hat u_i)\dot u_i^*,
       \eeq
       where $\eta''(\hat u_i)$ is the entropy Hessian evaluated
       at a convex combination $\hat u_i$ of $u_i$ and $u_i^{\mathrm{SSP}}$.

       By virtue of \eqref{semi-ineq}, the semi-discrete entropy fix
       guarantees the validity of
       $$
m_iv_i^T\dot u_i^*\le \sum_{j\in\mathcal N_i\backslash\{i\}}[G_{ij}-
   (\mathbf{q}_j-\mathbf{q}_i)\cdot\bc_{ij}]
       $$
under the sufficient condition (sum of inequalities \eqref{estab} over $j\in\mathcal N_i
\backslash\{i\}$)
\beq\label{ES1}
P_i^{\mathrm{SD}}\le
    \sum_{j\in\mathcal N_i\backslash\{i\}} (\boldsymbol{\psi}_j
  -\boldsymbol{\psi}_i) \cdot\mathbf{c}_{ij},
\eeq     
where
$$
P_i^{\mathrm{SD}}=\sum_{j\in\mathcal N_i\backslash\{i\}}
     \frac{(v_i-v_j)^T}2[d_{ij}(u_j-u_i)-
       (\mathbf{f}_j+\mathbf{f}_i) \cdot\mathbf{c}_{ij}
       +\alpha_{ij}f_{ij}^*]
$$
is the rate of entropy production by the spatial semi-discretization.
The second term on the right-hand side of \eqref{ent-taylor}
 is the rate of extra entropy production by the forward
 Euler time discretization.
           The
       fully discrete scheme is at least entropy conservative if
       \beq\label{ES2}
P_i^{\mathrm{SD}}+
       \frac{\Delta t}{2}m_i\,
       (\dot u_i^*)^T\eta''(\hat u_i)\dot u_i^*
   \le   P_i^{\mathrm{FD}}
\le  \sum_{j\in\mathcal N_i\backslash\{i\}}[(\boldsymbol{\psi}_j
  -\boldsymbol{\psi}_i) \cdot\mathbf{c}_{ij}+\dot G_{ij}]
\eeq
for an entropy production term $P_i^{\mathrm{FD}}$
and a flux $\dot G_{ij}$ to be defined below.
Under this condition, which may be more restrictive than
\eqref{ES1}, the fully discrete entropy inequality
 \beq\label{ent-fd}
 \eta(u_i^{\mathrm{SSP}})\le \eta(u_i)+\frac{\Delta t}{m_i}
    \sum_{j\in\mathcal N_i\backslash\{i\}}[G_{ij}+\dot G_{ij}-
   (\mathbf{q}_j-\mathbf{q}_i)\cdot\bc_{ij}]
       \eeq
 follows from \eqref{semi-ineq} and \eqref{ent-taylor}.
 Given an array of fluxes $f_{ij}^*$ preconstrained to satisfy
 \eqref{estab}, condition \eqref{ES2} can be enforced using
 correction factors $\alpha_{ij}\in[0,1]$. To avoid
 dependence of $\alpha_{ij}$ on the unknown state $\hat u$, we define
 the entropy production bound
$$
       P_i^{\mathrm{FD}}:=P_i^{\mathrm{SD}}+
       \frac{\Delta t}{2}m_i\,
       \langle\dot u_i^*,\dot u_i^*\rangle_{\eta''_i}
$$
using a generic
       Hessian-induced scalar product $\langle\cdot,\cdot\rangle_{\eta''_i}$.
       For a scalar conservation law, a~natural choice is
   $\langle f,g\rangle_{\eta_i''}:=\eta''_{\max}fg$,  where
       $\eta''_{\max}=\max\{\eta''(\hat u)\,:\,
       \hat u\in \mathcal G\}>0$. The
       scalar product of the bound $P_i^{\mathrm{FD}}$ for a hyperbolic system
       could formally be defined using the Hessian $\eta''(\hat u)$ evaluated
       at a state $\hat u\in\mathcal G$ such that $\dot u^T\eta''(\hat u)
       \dot u\ge \dot u^T\eta''(u)\dot u$ for all $u\in\mathcal G$
       and arbitrary $\dot u$. The existence of such a
       maximizer was
       conjectured in \cite{merriam1989}. In this work, we
approximate the unknown state $\hat u_i$ by $u_i$ and 
use $\langle f,g\rangle_{\eta_i''}:=f^T\eta''(u_i)g$ for systems.

Let the flux $\dot G_{ij}$ of the fully discrete entropy stability
condition \eqref{ES2} be defined by
$$
\dot G_{ij}=\frac{\Delta t}2\langle\dot u_i^L+\dot u_j^L,
\alpha_{ij}f_{ij}^*\rangle_{\eta''_i}.
$$
To transform \eqref{ES2} into a condition for calculating
$\alpha_{ij}$, we consider the decomposition
       $$
       \dot u_i^*=\dot u_i^L+\frac{1}{m_i}\sum_{j\in\mathcal N_i\backslash\{i\}}
       \alpha_{ij}f_{ij}^*,
       $$
       where
       \beq\label{udotL}
\dot u_i^L=\frac{1}{m_i}\sum_{j\in\mathcal N_i\backslash\{i\}}
     [d_{ij}(u_j-u_i)-
       (\mathbf{f}_j-\mathbf{f}_i) \cdot\mathbf{c}_{ij}]
     \eeq
          is the low-order approximation corresponding to 
$\alpha_{ij}=0\ \forall j\in\mathcal N_i\backslash\{i\}$. We have
     \begin{align*}
     \langle\dot u_i^*,\dot u_i^*\rangle_{\eta''_i}&=
       \langle\dot u_i^L,\dot u_i^L\rangle_{\eta''_i}+  
2\langle\dot u_i^L,\dot u_i^*-\dot u_i^L\rangle_{\eta''_i}
+\langle\dot u_i^*-\dot u_i^L,\dot u_i^*-\dot u_i^L\rangle_{\eta''_i}
=\langle\dot u_i^L,\dot u_i^L\rangle_{\eta''_i}\\&+
\frac{1}{m_i}\sum_{j\in\mathcal N_i\backslash\{i\}}
2\langle\dot u_i^L,\alpha_{ij}f_{ij}^*\rangle_{\eta''_i}+\frac{1}{m_i^2}
\left\langle
\sum_{j\in\mathcal N_i\backslash\{i\}}\alpha_{ij}f_{ij}^*,
\sum_{k\in\mathcal N_i\backslash\{i\}}
\alpha_{ik}f_{ik}^*\right\rangle_{\eta''_i}\\
&=\langle\dot u_i^L,\dot u_i^L\rangle_{\eta''_i}+
\frac{1}{m_i}\sum_{j\in\mathcal N_i\backslash\{i\}}
[ \langle\dot u_i^L+\dot u_j^L,\alpha_{ij}f_{ij}^*\rangle_{\eta''_i}
  +  \langle\dot u_i^L-\dot u_j^L,\alpha_{ij}f_{ij}^*\rangle_{\eta''_i}]\\
&+\frac{1}{m_i^2}
\left\langle
\sum_{j\in\mathcal N_i\backslash\{i\}}\alpha_{ij}f_{ij}^*,
\sum_{k\in\mathcal N_i\backslash\{i\}}
\alpha_{ik}f_{ik}^*\right\rangle_{\eta''_i}.
     \end{align*}
 Invoking the definition of $\dot G_{ij}$, the
 stability condition \eqref{ES2} can be written as
 \begin{align}
 \sum_{j\in\mathcal N_i\backslash\{i\}}\alpha_{ij}\left[
\frac{(v_i-v_j)^T}2f_{ij}^*\right. &+\left.\frac{\Delta t}2
\langle\dot u_i^L-\dot u_j^L,f_{ij}^*\rangle_{\eta''_i} \right]\nonumber\\
   &+\frac{\Delta t}{2m_i}
\left\langle
\sum_{j\in\mathcal N_i\backslash\{i\}}\alpha_{ij}f_{ij}^*,
\sum_{k\in\mathcal N_i\backslash\{i\}}
\alpha_{ik}f_{ik}^*\right\rangle_{\eta''_i}\nonumber\\
&\le \sum_{j\in\mathcal N_i\backslash\{i\}}Q_{ij}^{\mathrm{EC}}
-\frac{\Delta t}2m_i\langle\dot u_i^L,\dot u_i^L\rangle_{\eta''_i},\label{ES3}
\end{align}
where $Q_{ij}^{\mathrm{EC}}$ is the entropy-conservative bound
defined by \eqref{Qij-EC}. Note that if all terms
proportional to $\Delta t$ are set to zero, \eqref{ES3}
reduces to the sum of \eqref{estab} over
$j\in\mathcal N_i\backslash\{i\}$.

      The triangle inequality for the norm induced by the
      scalar product $\langle\cdot,\cdot\rangle_{\eta''_i}$
      yields 
  $$
\left\langle
\sum_{j\in\mathcal N_i\backslash\{i\}}\alpha_{ij}f_{ij}^*,
\sum_{k\in\mathcal N_i\backslash\{i\}}
\alpha_{ik}f_{ik}^*\right\rangle_{\eta''_i}\le\left(
\sum_{j\in\mathcal N_i\backslash\{i\}}\alpha_{ij}\sqrt{\langle
f_{ij}^*,f_{ij}^*\rangle_{\eta''_i}}\right)^2.
$$
It follows that a sufficient condition for the
  validity of \eqref{ES3} is given by
\begin{align}
    \sum_{j\in\mathcal N_i\backslash\{i\}}\alpha_{ij}
  \left[\frac{(v_i-v_j)^T}2f_{ij}^*\right.&+\left.\frac{\Delta t}{2}
    \langle \dot u_i^L-\dot u_j^L,f_{ij}^*\rangle_{\eta''_i}
   \right] +\frac{\Delta t}{2m_i}\left(
\sum_{j\in\mathcal N_i\backslash\{i\}}\alpha_{ij}\sqrt{\langle
  f_{ij}^*,f_{ij}^*\rangle_{\eta''_i}}\right)^2\nonumber\\
&\le \sum_{j\in\mathcal N_i\backslash\{i\}}Q_{ij}
-\frac{\Delta t}2m_i\langle\dot u_i^L,\dot u_i^L\rangle_{\eta''_i}.\label{ES4}
\end{align}
It can be configured to use the entropy-conservative
bound $Q_{ij}=Q_{ij}^{\mathrm{EC}}$, as defined by \eqref{Qij-EC}, or the
entropy-dissipative bound $Q_{ij}=Q_{ij}^{\mathrm{ED}}\le Q_{ij}^{\mathrm{EC}}$,
as defined by \eqref{Qij-ED}.

Following the derivation of the pressure limiter proposed in
\cite{lohmann2016}, we estimate the left-hand side of \eqref{ES4} 
using the fact that $\alpha_{ij}^2\le\alpha_{ij}$ for
$\alpha_{ij}\in[0,1]$. The correction factors $\alpha_{ij}$
of the fully discrete entropy fix
can then be calculated using the auxiliary quantities
$$
P_i=\sum_{j\in\mathcal N_i\backslash\{i\}}
\max\left\{0,
\frac{(v_i-v_j)^T}2f_{ij}^*+\frac{\Delta t}{2}
    \langle \dot u_i^L-\dot u_j^L,f_{ij}^*\rangle_{\eta''_i}
\right\} +\frac{\Delta t}{2m_i}\left(
\sum_{j\in\mathcal N_i\backslash\{i\}}\sqrt{\langle
  f_{ij}^*,f_{ij}^*\rangle_{\eta''_i}}\right)^2,
$$
$$
Q_i=\sum_{j\in\mathcal N_i\backslash\{i\}}Q_{ij}
-\frac{\Delta t}2m_i\langle\dot u_i^L,\dot u_i^L\rangle_{\eta''_i},\qquad
R_i=\min\left\{1,\frac{\max\{0,Q_i\}}{P_i}\right\}.$$
It is easy to verify that \eqref{ES4} holds if
$\alpha_{ij}\le R_i$ for all $j\in\mathcal N_i\backslash\{i\}$.
The choice of $\alpha_{ij}=\alpha_{ji}$ should also guarantee
the validity of \eqref{ES4} for node $j$.
Thus we limit $f_{ij}^*=-f_{ji}^*$ using
\beq\label{alpha-FD}
\alpha_{ij}=\min\{R_i,R_j\}.
\eeq
By the Taylor theorem, 
there is an intermediate state $\hat u$ such that
$(v_i-v_j)^T(u_j-u_i)=(u_j-u_i)^T\eta''(\hat u)(u_i-u_j)<0$
for $u_i\ne u_j$ and any convex entropy $\eta(u)$.
It follows that $Q_i\ge 0$ for $Q_{ij}$ defined by
\eqref{Qij-EC} or \eqref{Qij-ED}, provided that the
artificial viscosity coefficients $d_{ij}$ are chosen
sufficiently large. To avoid a priori verification for
$d_{ij}$ defined by \eqref{gms}, we use $\max\{0,Q_i\}$
to calculate $R_i$. In the case $Q_i\le 0$, formula
\eqref{alpha-FD} produces $\alpha_{ij}=0$, and a
fully discrete entropy inequality 
follows from the analysis of the ALF method in \cite{guermond2016}.

\begin{rmk}
  Since we assumed that the fluxes $f_{ij}^*$ are prelimited to
  satisfy \eqref{estab}, i.e., that the semi-discrete entropy
  fix is performed prior to the fully discrete one, the
  optional use of the latter
  can only make the AFC scheme more entropy dissipative.
\end{rmk}
  
\begin{rmk}\label{rmk:fde}
Condition \eqref{ES3} ensures fully discrete entropy stability
of each SSP Runge--Kutta stage. The rate of temporal entropy
production is proportional to the time step $\Delta t$, so the
fully discrete entropy fix may cause a loss of second-order
accuracy in applications to problems with smooth solutions.
The resolution of shocks is not affected because the spatial
discretization error dominates in the presence of shocks.
Lozano \cite{lozano2018} showed that high-order explicit
Runge--Kutta methods produce much less entropy than individual
stages. However, there is no simple way to limit the rate of
entropy production in a manner that exploits possible cancellation
effects. Moreover, the cost of a sophisticated second-order
explicit fix may be higher than that of the
iterative correction procedure proposed in Section \ref{sec:imp}.
  \end{rmk}

\subsection{Fully discrete implicit correction}
\label{sec:imp}

Let us now discretize \eqref{semi-afc-fix} in time using
a general $S$-stage Runge--Kutta method and constrain the final
stage in an iterative manner using the representation
\beq\label{full-rk}
u_i^{n+1}=u_i^n+\frac{\Delta t}{m_i} \sum_{j\in\mathcal N_i\backslash\{i\}}
[d_{ij}^{n+1}(u_j^{n+1}-u_i^{n+1})-(\mathbf{f}(u^{n+1}_j)-\mathbf{f}(u^{n+1}_i))\cdot
  \bc_{ij}+\alpha_{ij}^{n+1}f_{ij}^{n+1}],
\eeq
where $\alpha_{ij}^{n+1}\in[0,1]$ are implicitly defined
correction factors.
The space-time target flux
\begin{align}
  f_{ij}^{n+1}&=d_{ij}^{n+1}(u_i^{n+1}-u_j^{n+1})+
  (\mathbf{f}(u^{n+1}_j)+\mathbf{f}(u^{n+1}_i))\cdot
  \bc_{ij}\nonumber\\
  &+\sum_{s=1}^Sb_s[d_{ij}^{(s)}(u_j^{(s)}-u_i^{(s)})-
    (\mathbf{f}(u^{(s)}_j)+\mathbf{f}(u^{(s)}_i))\cdot
    \bc_{ij}+\alpha_{ij}^{(s)}
    f_{ij}^*(u^{(s)})]\label{flux-RK}    
\end{align}
is defined using the Butcher weights $b_s$ and property
\eqref{props-cij}. The intermediate stage approximations
$u^{(s)}$ are calculated using the bound-preserving MCL
limiter. The use of the semi-discrete entropy fix at
intermediate stages is optional. The high-order Runge--Kutta
time discretization of \eqref{semi-afc-fix} is recovered
in the case $\alpha_{ij}^{n+1}=1\ \forall j\in
\mathcal N_i\backslash\{i\}$.

In the process of flux correction for the final stage
\eqref{full-rk}, the IDP property is enforced
(as in \cite{kuzmin-arxiv}) using the
MCL limiter for $f_{ij}^{n+1}$. To prevent a loss of
accuracy at this stage, the bounds 
should be global (as in \cite{kuzmin-arxiv}) or defined using
all states $u_j^n$ that may influence $u_i^{n+1}$ if
time integration is performed using the high-order
Runge--Kutta scheme.

The entropy correction factor $\alpha_{ij}^{n+1}$ to be used
in \eqref{full-rk} is defined by \eqref{alpha}. This
definition ensures fully discrete entropy stability without
additional fixes since \eqref{full-rk} has the structure of the
backward Euler (BE) method for an AFC scheme of the
form \eqref{semi-afc-fix}.
Note that the first-order
BE scheme for \eqref{semi-afc-fix} would
use $f_{ij}^*(u^{n+1})$ instead of
$f_{ij}^{n+1}$ defined by \eqref{flux-RK}.

In the numerical studies of Section \ref{sec:results}, we use the
second-order SSP-RK method 
$$
u^{(1)}=u^n,\qquad
u^{(2)}=u^n+\Delta t\dot u^{(1)},
\qquad u^{n+1}=u^n+\Delta t\frac{\dot u^{(1)}+\dot u^{(2)}}2
\qquad 
$$
and solve \eqref{full-rk} using the simple fixed-point
iteration $u_i^{[k+1]}=u_i^n+\Delta t \dot u_i^{[k]},$ where
$$
 \dot u_i^{[k]}=\frac{1}{m_i} \sum_{j\in\mathcal N_i\backslash\{i\}}
      [d_{ij}^{[k]}(u_j^{[k]}-u_i^{[k]})-(\mathbf{f}(u^{[k]}_j)-
          \mathbf{f}(u^{[k]}_i))\cdot
  \bc_{ij}+\alpha_{ij}^{[k]}f_{ij}^{[k]}].
$$
Fully discrete entropy stability of the converged
solution $u_i^{n+1}$ follows from
 \beq
    \frac{\eta(u_i^{n+1})-\eta(u_i^n)}{\Delta t}
    =(v_i^{n+1})^T\dot u_i^{n+1}
       - \frac{\Delta t}{2}
      (\dot u_i^{n+1})^T\eta''(\hat u_i)\dot u_i^{n+1},
       \eeq
       where $\dot u_i^{n+1}$ satisfies a semi-discrete entropy inequality
       (cf. \cite{merriam1989}, Section 4.3). The IDP property is
       guaranteed for any $\Delta t$, even if the
       underlying RK scheme is explicit. However, instability of the
       baseline discretization may trigger aggressive flux limiting.
       Therefore, the RK method corresponding to $\alpha_{ij}^{n+1}=1\
       \forall j\in\mathcal N_i\backslash\{i\}$
       should be at least linearly stable. Moreover, the
       time step should be small enough for the
       fixed-point iteration to converge.

\section{Convergence to entropy solutions}
\label{sec:lw}

To further motivate the use of bound-preserving and entropy-stable schemes, we prove a finite element version of the Lax--Wendroff theorem in this section. To that end, we introduce some additional notation and auxiliary results that are needed in the proof.

First, we recognize that approximations $u_h$ produced by the aforementioned fully discrete schemes have only been defined at discrete time instants\ie for $t \in \{0, \Delta t, 2 \Delta t, \ldots\}$ so far. However, for analysis purposes, we interpret the numerical solution $u_h$ as an element of $C^0(\Omega \times \R^+_0)$. Using linear interpolation between the time levels, we set
\begin{equation}\label{EQ:space-time}
 u_h(\vec x, t) =  \frac{ ((n+1) \Delta t - t) u_h(\vec x, n \Delta t) + (t - n \Delta t) u_h(\vec x, (n+1) \Delta t)}{\Delta t}
\end{equation}
for $\vec x \in \bar\Omega$ and $t \in (n \Delta t, (n+1) \Delta t)$. This extension enables us to interpret $u_h(\cdot,t)$ as a piecewise-linear or multilinear finite element function on $\mathcal T_h$ at time $t\in\R_0^+$ and $u_h(\mathbf{x},\cdot)$ as a continuous piecewise-linear interpolant of the discrete values $u_h(\mathbf{x},n\Delta t),\ n=0,1,2,\ldots$ at the nodes of a uniform time grid (denoted by $\timeMesh$) for $\R^0_+$ at $\mathbf{x}\in \bar\Omega$. In this sense, we have extended $u_h$ to be a continuous space-time finite element approximation.

Second, $u_h$ depends both on the mesh size $h$ and on the time step $\Delta t$. The sequence of approximations $u_h$ should converge to a weak solution $u$ of \eqref{ibvp} in the limit $h \to 0$ and $\Delta t \to 0$. In view of the CFL condition, we choose sequences $(h_k)$ and $((\Delta t)_k)$ such that
\begin{equation}\label{EQ:sequence}
  h_k \to 0\quad \text{ as }  k \to \infty\quad \implies\quad
   (\Delta t)_k \to 0 \quad \text{ as } \ k \to \infty.
\end{equation}
This refinement strategy makes it possible to avoid 
indeterminancies
in ratios involving $h_k$ and $(\Delta t)_k$.
Moreover, it is
consistent with the common practice of choosing a time step $(\Delta t)_k$ that
depends on the mesh size $h_k$ in real-world simulations. 

Let us now define some functional spaces and discuss their properties. In this section, $\tilde V_h$ denotes the finite element space of continuous piecewise-linear functions for a conforming triangulation $\mathcal T_h$ of a domain $\Omega\subset\R^d$. As before, we assume that $\Omega$ has periodic boundaries. The spatial interpolation operator $\tilde \interpol_h: C^0(\Omega) \to \tilde V_h$ is defined in the usual sense. The Bochner space $V_{h_k} := C^0_c(\R^+_0;\tilde V_{h_k})$ consists of functions that are continuous and have compact support as functions of $t\in\R^+_0$, while being continuous finite element functions of the space variable $\mathbf{x}\in\bar\Omega$. The interpolation operator for $V_{h_k}$ is defined by
\begin{align}
 \interpol_{h_k}: C^2_c (\Omega \times \R^+_0)  \to V_{h_k},\qquad
 \varphi(\cdot, t)  \mapsto \tilde \interpol_{h_k} \varphi(\cdot, t) \quad \forall t \in \R^+_0.
\end{align}
That is, for any fixed time $t \in \R^+_0$, the operator $\interpol_{h_k}$ interpolates a given function to the finite element space $\tilde V_h$. The subscript $c$ is again used to indicate that functions belonging to the space have compact support in time. This has several immediate consequences:
\begin{enumerate}
\item Interpolation and time derivatives commute:
\smash{\mbox{  $\tfrac{\partial}{\partial t} \interpol_{h_k} \varphi = \interpol_{h_k} \tfrac{\partial\varphi}{\partial t} $, $\tfrac{\partial^2}{\partial t^2} \interpol_{h_k} \varphi = \interpol_{h_k} \tfrac{\partial^2\varphi}{\partial t^2} $.}} This property also implies that $\interpol_{h_k}$ is well defined.
 \item We have \smash{\mbox{$\| \interpol_{h_k} \varphi(\cdot, t) \|_{L^\infty(\Omega)} \le \| \varphi(\cdot, t) \|_{L^\infty(\Omega)}$, $\|\tfrac{\partial}{\partial t} \interpol_{h_k} \varphi(\cdot, t) \|_{L^\infty(\Omega)} \le \| \tfrac{\partial\varphi}{\partial t} (\cdot, t) \|_{L^\infty(\Omega)}$}}, and \smash{\mbox{$\| \tfrac{\partial^2}{\partial t^2} \interpol_{h_k} \varphi(\cdot, t) \|_{L^\infty(\Omega)} \le \| \tfrac{\partial^2\varphi}{\partial t^2} (\cdot, t) \|_{L^\infty(\Omega)}$}} for all $t \in \R^+$.
 \item There exists a constant $C_\varphi > 0$ such that 
 \begin{align*}
  \| \interpol_{h_k} \varphi \|_{L^\infty(\Omega\times \R^+)} + \| \varphi \|_{L^\infty(\Omega \times \R^+)} &+ \| \tfrac{\partial}{\partial t} \interpol_{h_k} \varphi \|_{L^\infty(\Omega\times \R^+)} + \| \tfrac{\partial\varphi}{\partial t}  \|_{L^\infty(\Omega \times \R^+)} \\
  &+ \| \tfrac{\partial^2\varphi}{\partial t^2} \interpol_{h_k}  \|_{L^\infty(\Omega\times \R^+)} + \| \tfrac{\partial^2\varphi}{\partial t^2}  \|_{L^\infty(\Omega \times \R^+)}  \le C_\varphi.
  \end{align*}
\end{enumerate}
Finally, let us define the total variation of a space-time finite element function $u_{h_k}$ as follows:
\begin{equation}
 \operatorname{TV}(u_{h_k}) := \sup_{t \in \R^+_0} \sum_{K \in \mathcal T_h} \left[\sup \left\{ u_{h_k}(t, \vec x_1) - u_{h_k}(t,\vec x_2) \mid \vec x_1, \vec x_2 \in K \right\}\right]^2.
\end{equation}
This definition generalizes the one given in \cite[(12.40)]{leveque1992} to functions from the space $V_{h_k}$.
\smallskip

With these preliminaries, we are now ready to prove a finite element version
of the Lax--Wendroff theorem \cite[Theorem 12.1]{leveque1992} 
about convergence to weak solutions.
\begin{theorem}[Convergence of finite element schemes to weak solutions]\label{TH:LW}
 Consider\\ problem \eqref{ibvp} with periodic boundary conditions. Suppose that $u_0 \in C^0(\bar\Omega)$ and\\ $\vec f \in C^{0}(\R)^d$ is Lipschitz with constant $C_\vec f > 0$. Initialize the numerical solutions $u_{h_k}$\\ by $u_{h_k}(\cdot, 0) = \interpol_{h_k} u_0$ and evolve them using a fully discrete scheme of the form
 \begin{align}
   \sum_{i=1}^{N_{h_k}}(u^{n+1}_{h_k,i} - u^{n}_{h_k,i}) 
   \int_\Omega \tilde \phi_{k,i} \dx&+
 (\Delta t)_k \int_\Omega\tilde \phi_{h_k} (\nabla \cdot \vec f^n_{h_k}) \dx \nonumber\\
   &=
   (\Delta t)_k\mathbb S_k (u_{h_k}^n, \tilde \phi_{h_k})
   \qquad \forall \tilde \varphi_{h_k} \in \tilde V_{h_k},\ n \in \N_0,\label{EQ:scheme}
 \end{align}
 where $\mathbb S_k$ is a stabilization term.
 Assume that there exist a time $T > 0$, a~function\\ $u \in L^2(\Omega \times (0,T))$,  and a constant $C_u>0$ independent of $k\in\N$ such that
 \begin{gather}\label{EQ:u_conv}
  \| u_{h_k} - u \|_{L^2(\Omega \times (0,T))} \to 0 \qquad \text{as}\ k \to \infty,\\
  \| u_{h_k} \|_{L^\infty(\Omega \times (0,T))} + \operatorname{TV}(u_{h_k}) \le C_u.
  \end{gather}
  Furthermore, assume that
   $\sum_{i=1}^{N_{h_k}}\mathbb S_{k,i}=0$ and the Ritz projections
\beq\label{ritz}
  \int_\Omega\nabla \phi_{h_k}\cdot \nabla q_{h_k}^n\dx=
\mathbb S_k(u^n_k,\phi_{h_k})\qquad \forall \phi_{h_k}\in \tilde V_{h_k}
\eeq
produce flux potentials $q_{h_k}^n$ such that
 \begin{equation}
 \| \nabla q_{h_k}^n\|_{L^2(\Omega)} \to 0 \qquad \text{ as } k \to \infty.\label{ritz-est}
 \end{equation}
 Then $u$ is a weak solution to \eqref{ibvp} in the sense that
 \begin{equation}\label{EQ:analytic}
  \int_0^T \int_\Omega \left[\tfrac{\partial \phi}{\partial t} u + \nabla \phi\cdot \vec f(u)\right] \dx \dt = - \int_\Omega \phi(\vec x, 0) u(\vec x, 0) \dx
 \end{equation}
 for all test functions
 $\phi \in C^2_c(\Omega \times [0,T))$ (compact support in time, periodic in space).
\end{theorem}
\begin{rmk}
 The left-hand side of \eqref{EQ:scheme} corresponds to a fully discrete version of the lumped-mass
 group finite element method \eqref{semi-group}. Inexact quadrature
introduces a consistency error which vanishes as $k \to \infty$.
The global conservation property $\sum_{i=1}^{N_{h_k}}\mathbb S_{k,i}=0$
ensures solvability of \eqref{ritz}. The flux potential $q_{h_k}^n$ is defined up to a constant which
has no influence on the value of $\nabla q_{h_k}^n$. The
use of $\mathbb S_k$ in \eqref{EQ:scheme} has the same effect as
addition of $-\nabla q_{h_k}^n$ to $\vec f_{h_k}^n$.
\end{rmk}
\begin{proof}[Proof of Theorem \ref{TH:LW}]
  The test functions $\tilde \varphi_{h_k}$ of the discrete problem \eqref{EQ:scheme} are arbitrary elements of $\tilde V_{h_k}$. In particular, they may represent projections of $\phi \in C^2_c(\Omega \times [0,T))$ into $\tilde V_{h_k}$ at discrete time levels $n(\Delta t)_k$ or intermediate time instants $t \in (n (\Delta t)_k, (n+1) (\Delta t)_k)$. Let us construct a sequence of functions $\varphi_{h_k} \in V_{h_k}$ that converges to $\varphi$. Then we may choose $\varphi_{h_k}^n := \varphi_{h_k}(\cdot,n(\Delta t)_k)$ to be the test function $\tilde \varphi_{h_k}\in\tilde V_{h_k}$ for \eqref{EQ:scheme}. 
    

 We need to show that \eqref{EQ:analytic} holds for functions $u$ and $\phi$, which are supposed to be limits of the sequences $(u_{h_k})$ and $(\phi_{h_k})$. The space-time finite element approximations $u_{h_k}$ are defined by the numerical scheme, while $\phi_{h_k} \in V_{h_k}$ can be chosen arbitrarily. Using the test functions $\phi_{h_k} := \interpol_{h_k} \varphi$, we will show that $\phi_{h_k} \to \phi$ in an appropriate sense.

 Let us first cast \eqref{EQ:scheme} into a form that better resembles \eqref{EQ:analytic}. Summing the discretized equations over all time steps and using transformations to be explained below, we obtain
 \begin{align}\begin{split}\label{EQ:conv}
0 = & - \underbrace{ \int_\Omega u^0_{h_k} \phi_{h_k}^0 \dx }_{ =: \Xi^k_1 } - \underbrace{ (\Delta t)_k \sum_{n=1}^\infty  \int_\Omega \frac{\phi_{h_k}^{n} - \phi_{h_k}^{n-1}}{(\Delta t)_k} u^n_{h_k} \dx }_{ =: \Xi_2^k}\\ 
&- \underbrace{ (\Delta t)_k \sum_{n=0}^\infty \int_\Omega
  \nabla \phi^n_{h_k} \cdot (\vec f_{h_k}-\nabla q_{h_k})
  \dx }_{ =: \Xi_3^k }\\ &+
   \underbrace{\sum_{n=0}^\infty \left[ \sum_{i=1}^{N_{h_k}}(u^{n+1}_{h_k,i} - u^{n}_{h_k,i}) 
   \int_\Omega \phi^n_{h_k,i} \dx- \int_\Omega (u_{h_k}^{n+1} - u_{h_k}^n)\varphi^n_{h_k}\dx \right].
  }_{ =: \Xi_4^k.}
 \end{split}\end{align}
 To derive the term $\Xi^k_3$, we used integration by parts in the volume integral involving the divergence of $\vec f^n_{h_k}$, added the stabilization term $\mathbb{S}_k$, and expressed the result in terms of flux potentials $q_{h_k}$ using \eqref{ritz}. The terms $\Xi^k_1$ and $\Xi_2^k$ were obtained using summation by parts, a~discrete version of integration by parts for sums with a finite number of nonzero terms. For a more detailed description of this procedure, we refer to \cite[(12.45)--(12.47)]{leveque1992}.

 The remainder of the proof consists of four steps. In each step, we show that a~term present in \eqref{EQ:conv} converges to its counterpart in \eqref{EQ:analytic} or to zero in the case of $\Xi_4^k$.

 \paragraph{Step 1: $\Xi_1^k \to \int_\Omega \phi(\vec x, 0) u(\vec x, 0) \dx$} This follows directly from  $u_{h_k}(\cdot,0) = \interpol_{h_k} u_0$, $\phi_{h_k}(\cdot,0) = \interpol_{h_k} \phi(\cdot,0)$ and standard convergence results for linear interpolation operators.

 \paragraph{Step 2: $\Xi_2^k \to \int_0^T \int_\Omega \tfrac{\partial\phi}{\partial t}  u \dx \dt$} We observe that $\int_0^T \int_\Omega\tfrac{\partial \phi}{\partial t} u \dx \dt - \Xi^2_k$ can be written as
 \begin{subequations}\label{EQ:xi_2}\begin{align}
  \int_0^T \int_\Omega & \frac{\partial \varphi }{\partial t} u \dx \dt - \int_0^T \int_{\Omega} \frac{\partial \phi_{h_k}}{\partial t}  u_{h_k} \dx \dt\label{EQ:xi_21}\\
  & + \int_0^T \int_{\Omega} \frac{\partial \phi_{h_k}}{\partial t} u_{h_k} \dx \dt - (\Delta t)_k \sum_{n=1}^\infty \int_{\Omega} \frac{\partial \phi^n_{h_k} }{\partial t}u^n_{h_k} \, \dx \dt \label{EQ:xi_22}\\ 
  & + (\Delta t)_k \sum_{n=1}^\infty \int_\Omega \frac{\partial\phi^n_{h_k} }{\partial t} u^n_{h_k} \, \dx - (\Delta t)_k \sum_{n=1}^\infty \int_{\Omega} \frac{\phi_{h_k}^{n} - \phi_{h_k}^{n-1}}{(\Delta t)_k} u^n_{h_k} \, \dx.\label{EQ:xi_23}
 \end{align}\end{subequations}
 To prove the desired result, we show that the three lines of \eqref{EQ:xi_2} go to zero:
 \begin{enumerate}[a)]
 \item We know that $u_{h_k} \to u$ in the sense of \eqref{EQ:u_conv}, and that $\tfrac{\partial\phi_{h_k}}{\partial t}  \to \tfrac{\partial\varphi}{\partial t}$
by the approximation property of the interpolation operator
that we used to construct $\phi_{h_k}$.
\item This line is the error of the trapezoidal quadrature rule.
  We estimate it as follows:
  \begin{align*}
   | \eqref{EQ:xi_22} | \le & \| u_{h_k} \|_{L^\infty(\Omega \times \R^+)} \sum_{n=0}^\infty \int_{n(\Delta t)_k}^{(n+1)(\Delta t)_k} \int_{\Omega} \left| \tfrac{\partial\phi_{h_k}}{\partial t}  - \tfrac{\partial\phi^{n+1}_{h_k} }{\partial t} \right| \, \textup \dx \, \textup dt \\
   \le & C_u T |\Omega| \max_{n\in\N_0} \underbrace{ \left\| \tfrac{\partial\phi_{h_k}}{\partial t}  - \tfrac{\partial}{\partial t} \phi_{h_k}(\cdot,(n+1)(\Delta t)_k)  \right\|_{L^\infty(\Omega \times (n (\Delta t)_k, (n+1) (\Delta t)_k)} }_{ \le (\Delta t)_k \left\| \tfrac{\partial^2\varphi}{\partial t^2}  \right\|_{L^\infty(\Omega \times (n (\Delta t)_k, (n+1) (\Delta t)_k)} }.
  \end{align*}
 The first two estimates are consequences of H\"older's inequality. The underbraced inequality is a uniform bound for the remainder of a Taylor expansion with respect to time. The total quadrature error goes to zero as $(\Delta t)_k \to 0$.
\item A similar argument using a
 Taylor expansion with a remainder yields the estimate
  \begin{align*}
   | \eqref{EQ:xi_23} | \le & \left[\sum_{n=1}^{T / (\Delta t)_k} (\Delta t)_k \right] \| u_{h_k} \|_{L^\infty(\Omega \times \R^+)} \max_{n\in\N}  \int_{\Omega} \left| \tfrac{\partial\phi^n_{h_k}}{\partial t}  - \frac{\phi_{h_k}^{n} - \phi_{h_k}^{n-1}}{(\Delta t)_k} \right| \dx  \\
   \le & T C_u \max_{n\in\N} \int_{\Omega} \left| \tfrac{\partial\phi^n_{h_k} }{\partial t} - \tfrac{\partial \phi_{h_k}}{\partial t}(\vec x,\xi(\vec x)) \right| \dx,
  \end{align*}
  where $\xi(\vec x) \in (n(\Delta t)_k, (n+1) (\Delta t)_k)$. The last term goes to zero as shown above.
 \end{enumerate}

 \paragraph{Step 3: $\Xi^3_k \to \int_0^T \int_\Omega \nabla \phi\cdot \vec f(u) \dx \dt$} To
prove this, we write
 $\int_0^T \int_\Omega \nabla \varphi\cdot \vec f(u) \dx\ \dt - \Xi_3^k$ as
 \begin{subequations}\begin{align}
  \int_0^T \int_\Omega & \nabla \varphi\cdot \vec f(u) \dx \dt - \int_0^T \int_\Omega \nabla \phi \cdot \vec f(u_{h_k}) \dx \dt \label{EQ:s3_1}\\
  & + \int_0^T \int_\Omega \nabla \phi \cdot \vec f(u_{h_k}) \dx \dt - \int_0^T \int_\Omega \nabla \phi \cdot \vec f_{h_k} \dx \dt \\
  & + \int_0^T \int_\Omega \nabla \phi \cdot \vec f_{h_k} \dx \dt - \int_0^T \int_\Omega \nabla \phi \cdot (\vec f_{h_k} - \nabla q_{h_k} ) \dx \dt\\
  & + \int_0^T \int_\Omega \nabla \phi \cdot (\vec f_{h_k} - \nabla q_{h_k} ) \dx \dt - \int_0^T \int_\Omega \nabla \phi_{h_k} \cdot (\vec f_{h_k} - \nabla q_{h_k} ) \dx \dt\\
  & + \int_0^T \int_\Omega \nabla \phi_{h_k} \cdot (\vec f_{h_k} - \nabla q_{h_k} ) \dx \dt - (\Delta t)_k \sum_{n=0}^\infty \int_\Omega \nabla \phi^n_{h_k} \cdot (\vec f_{h_k} - \nabla q_{h_k} ) \dx,
 \end{align}\end{subequations}
 where $\mathcal N_K$ denotes the integer set in which the indices of nodes $\mathbf{x}_{j_k}$ belonging to $K\in \mathcal T_{h_k}$ are stored. Next, we use the following arguments to show that the five lines go to zero.
 \begin{enumerate}[a)]
 \item This follows immediately, since $\mathbf{f}$ is Lipschitz continuous and $u_{h_k} \to u$ in $L^2$ by assumption. It follows that $\| \vec f(u) - \vec f(u_{h_k}) \|_{L^2(\Omega \times (0,T))}\to 0$.
 
  \item This is the difference between the standard and group finite element formulations. It tends to zero for any $t \in (0,T)$ because
  \begin{align*}
    \sum_{K \in \mathcal T_{h_k}} \int_K &\left|\vec f(u_{h_k}) - \sum_{j\in\mathcal N_K} \vec f(u_{h_k}(\mathbf{x}_{j_k})) \phi_{kj} \right|^2 \dx 
   \le  \sum_{K \in \mathcal T_{h_k}}  |K| \max_{j\in\mathcal N_K} \| \vec f(u_{h_k}) - \vec f(u_{h_k}(\vec x_{j_k})) \|^2_{L^\infty(K)} \\
   &\qquad\le  \sum_{K \in \mathcal T_{h_k}} |K| C_{\vec f} \max_{j\in\mathcal N_K} \| u_{h_k} - u_{h_k}(\vec x_{j_k}) \|^2_{L^\infty(K)} 
   \le  \max_{K \in \mathcal T_{h_k}} |K| C_{\vec f} C_u.
  \end{align*}
  \item Using the Cauchy-Schwarz inequality and assumption \eqref{ritz-est}, we find that  the difference of the two integrals goes to zero.
  \item Since convergent sequences are bounded and  $\nabla \phi^n_{h_k} \to \nabla \phi$ in $L^2$, this line goes to zero too.
  \item The fact that this line goes to zero can be verified similarly to \eqref{EQ:xi_22}.
 \end{enumerate}
 
 \paragraph{Step 4: $\Xi_4^k \to 0$} This term goes to zero since the difference between the lumped and consistent mass version is the quadrature error of a low-order Newton--Cotes rule.
 \medskip
 
 The assertion of the Theorem follows from the convergence results for all steps.
\end{proof}
Theorem \ref{TH:LW} does not rule out convergence to a wrong weak solution. However,
if \eqref{EQ:scheme} is entropy stable \wrt an entropy pair, it can only
converge to a solution satisfying a continuous weak form of the corresponding
entropy inequality. This property of finite element approximations to hyperbolic
problems is guaranteed by the following theorem.

\begin{theorem}[Convergence of finite element schemes to entropy solutions]
  Under the assumptions of Theorem \ref{TH:LW}, let $\{\eta(u), \vec q(u)\}$ be an entropy pair such that $\eta \in C^{1}(\R)$, the corresponding entropy variable $v = \eta'(u)$ is Lipschitz, and so is $\vec q \in C^{0}(\Omega)^d$. Assume that
 \begin{align}
   \sum_{i=1}^{N_{h_k}} \left(\eta(u^{n+1}_{h_k,i}) - \eta(u^n_{h_k,i})\right) \int_\Omega \tilde \varphi_{h_k,i} \dx &+ (\Delta t)_k\int_\Omega \tilde \varphi_{h_k}(\nabla \cdot \vec q^n_{h_k}) \dx
   \nonumber\\
   &  \le (\Delta t)_k \mathbb S^\eta_k(\eta(u^n_{h_k}),\tilde \phi_{h_k}) \qquad
   \forall \tilde \varphi_k \in \tilde V_{h_k},\ n \in \N_0.\label{FDLW}
 \end{align}
 Furthermore, assume that $\sum_{k=1}^{N_{h_k}} \mathbb S^\eta_k=0$ and the Ritz projections
$$
 \int_\Omega\nabla \phi_{h_k}\cdot \nabla r_{h_k}^n\dx=
  \mathbb S^\eta_k(\eta(u^n_{h_k}),\phi_{h_k})\qquad
  \forall \phi_{h_k}\in \tilde V_{h_k}$$
 produce flux potentials $r_{h_k}^n$ such that
 \begin{equation*}
\| \nabla r_{h_k}^n\|_{L^2(\Omega)} \to 0 \qquad \text{ as } k \to \infty.
 \end{equation*}
 Then the weak entropy inequality
 \begin{equation}\label{EQ:weak_ent}
  \int_0^T \int_\Omega \left[\tfrac{\partial\varphi }{\partial t} \eta(u) + \nabla \varphi\cdot \vec q(u)\right] \dx \dt \ge - \int_\Omega \varphi(\vec x, 0) \eta(u(\vec x, 0)) \dx
 \end{equation}
 holds for all $\varphi \in C^2_c(\Omega \times [0,T); \R^+_0)$.
\end{theorem}
\begin{rmk}
  The weak form  \eqref{EQ:weak_ent} of the
  entropy inequality can be derived from \eqref{ent-ineq}
  using multiplication by a nonnegative, smooth test functions with compact support with respect to time, integration over the space-time domain $\Omega\times(0,T)$, and integration by parts. 
\end{rmk}
\begin{proof}
 The proof of this theorem is similar to that of Theorem \ref{TH:LW}. Thus, we skip the details of steps that involve the same arguments. 
 Summing the discrete entropy inequalities \eqref{FDLW} over
 all time steps and following the proof of Theorem \ref{TH:LW}, we arrive at 
 \begin{align*}
  0 \ge & - \underbrace{ \int_\Omega\eta(u^0_{h_k}) \phi_{h_k}^0 \dx }_{ =: \Xi^k_1 } - \underbrace{ (\Delta t)_k \sum_{n=1}^\infty  \int_\Omega \frac{\phi_{h_k}^{n} - \phi_{h_k}^{n-1}}{(\Delta t)_k} \eta(u^n_{h_k}) \dx }_{ =: \Xi_2^k} \\
  &- \underbrace{ (\Delta t)_k \sum_{n=0}^\infty \int_\Omega \nabla \phi^n_k
    \cdot(\vec q_{h_k}-\nabla r_{h_k})  \dx }_{ =: \Xi_3^k }\\
  & + \underbrace{\sum_{n=0}^\infty \left[ \sum_{i=1}^{N_{h_k}}(\eta(u^{n+1}_{h_k,i}) - \eta(u^{n}_{h_k,i})) \int_\Omega \phi^n_{h_k,i} \dx- \int_\Omega (\eta(u_{h_k}^{n+1}) - \eta(u_{h_k}^n))\varphi^n_{h_k}\dx \right].
  }_{ =: \Xi_4^k.}
 \end{align*}
 We choose $\phi_{h_k}$ as above and prove that the terms in the first line
 converge to their continuous counterparts. The convergence proofs for $\Xi_k^1$ and $\Xi^2_k$ repeat those in  Theorem \ref{TH:LW} using Lipschitz continuity to show
 that $\eta(u_{h_k}) \to \eta(u)$ in $L^2$ if $u_{h_k} \to u$ in $L^2$. The estimation of $\Xi_k^3$ is also performed along similar lines using $r_{h_k}$ and $\vec q(u_{h_k})$ in place of
 $q_{h_k}$ and $\vec f^n_k$. Again,  Lipschitz continuity implies that $\vec q(u_{h_k}) \to \vec q(u)$ in $L^2$ if $u_{h_k} \to u$ in $L^2$. The quadrature error
$\Xi_k^4$ converges to zero, which concludes the outline of the proof.
\end{proof}

\section{Numerical experiments}\label{sec:results}

To illustrate the numerical behavior of limiter-based entropy correction
procedures, we apply them to a suite of standard nonlinear test problems in this
section. Unless mentioned otherwise, we integrate in time using a
second-order explicit SSP Runge--Kutta scheme (Heun's method) and
perform algebraic flux correction using the target fluxes
\beq
f_{ij}=m_{ij}(\dot u_i^L-\dot u_j^L)+d_{ij}(u_i-u_j),\label{target}
\eeq
where $\dot u_i^L$ is defined by \eqref{udotL}. This choice
corresponds to a stabilized Galerkin approximation; see
\cite{kuzmin2020} for details. To better illustrate the
effectiveness of entropy fixes for systems, we also use Roe's
approximate Riemann solver (a generalization of which to
finite element discretizations can be found in \cite{kuzmin2012b,lyra2002,lohner2008})
as the target scheme in some examples.

In figures and descriptions of numerical results, we use various combinations of
the following acronyms to distinguish between different
methods under investigation:
\begin{itemize} 
\item HO: high-order scheme without flux correction;
\item LO: low-order ALF scheme defined by \eqref{semi-low};  
\item BP: bound-preserving AFC without entropy fixes;
\item EC: entropy fix using $Q_{ij}=Q_{ij}^{\mathrm{EC}}$ defined by \eqref{Qij-EC};
\item ED: entropy fix using $Q_{ij}=Q_{ij}^{\mathrm{ED}}$ defined by \eqref{Qij-ED};
\item SD: semi-discrete entropy fix of Section \ref{sec:semi};
\item FDE: fully discrete explicit fix of Section \ref{sec:exp};
\item FDI: fully discrete implicit fix of Section \ref{sec:imp}.
\end{itemize}
Additionally, the type of the target flux (GT:= Galerkin target,
RT:=Roe target) may be specified for a given AFC scheme. The default
is GT, i.e., $f_{ij}$ defined by \eqref{target}.

\subsection{One-dimensional KPP problem}

We begin with numerical studies for
one-dimensional scalar equations. The objective of the test
proposed in \cite{kurganov2007}
is to assess the ability of numerical methods to produce approximations that converge to the vanishing viscosity solution of the conservation law
\begin{align*}
\frac{\partial u}{\partial t} + \frac{\partial f(u)}{\partial x} = 0, \qquad f(u) = \begin{cases}
\frac1 4 u(1-u) & \mbox{if } u \le 0.5, \\
\frac 1 2 u(u-1) + \frac3{16} & \mbox{otherwise.}
\end{cases}
\end{align*}
Following Kurganov et al. \cite{kurganov2007}, we consider the following two Riemann problems:
\begin{align*}
\text{RP1:}\quad u_0(x) = \begin{cases}
0 & \mbox{if } x < 0.25, \\ 1 & \mbox{otherwise,} \end{cases} \qquad\qquad
\text{RP2:}\quad u_0(x) = \begin{cases}
1 & \mbox{if } x < 0.25, \\ 0 & \mbox{otherwise.} \end{cases}
\end{align*}
The corresponding vanishing viscosity solutions
\begin{align*}
\text{RP1:}\quad& u(x,t) = \begin{cases}
0 & \mbox{if } x < \frac{1 + (\sqrt{6}-2) t }4, \\
\frac 1 2 + \frac{x-1/4}{t} & \mbox{if } \frac{1 + (\sqrt{6}-2) t }4< x < (1+2t)/4, \\
1 & \mbox{otherwise,}
\end{cases}
\\
\text{RP2:}\quad& u(x,t) = \begin{cases}
1 & \mbox{if } x < \frac 1 4(1+(\sqrt 3 - 1)t), \\
\frac 1 2 - \frac{2(x-1/4)}{t} & \mbox{if } 
\frac 1 4(1+(\sqrt 3 - 1)t) < x < (2+t)/4, \\
0 & \mbox{otherwise}
\end{cases}
\end{align*}
can be derived using the family of Kruzkov entropy-entropy flux pairs $\eta_{\tilde u}(u) = |u-\tilde u|$, $q(u) = \mathrm{sgn}(u-\tilde u)(f(u)-f(\tilde u))$, where $\tilde u \in \R$ \cite{kurganov2007}.

\begin{figure}[ht!]
\centering
\begin{subfigure}[b]{0.49\textwidth}
\caption{RP1 at $t=1.0$, no fully discrete entropy fixes}
\includegraphics[width=\textwidth]{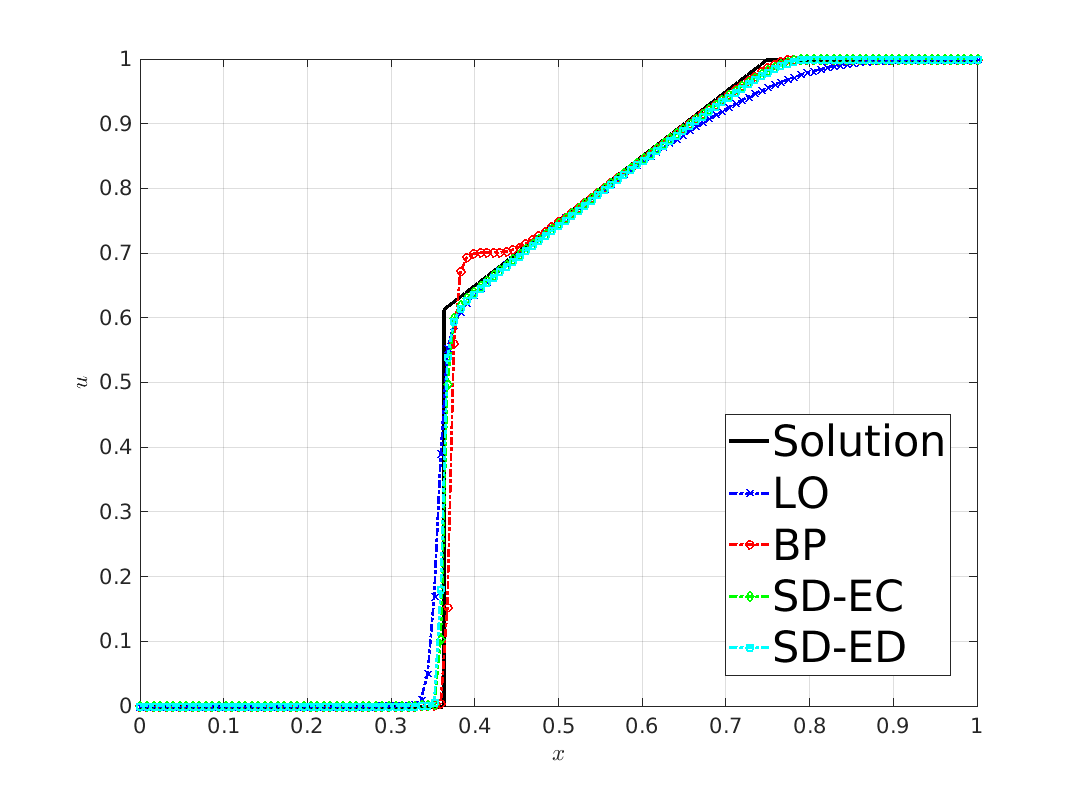}
\end{subfigure}
\begin{subfigure}[b]{0.49\textwidth}
\caption{RP1 at $t=1.0$ with fully discrete entropy fixes}
\includegraphics[width=\textwidth]{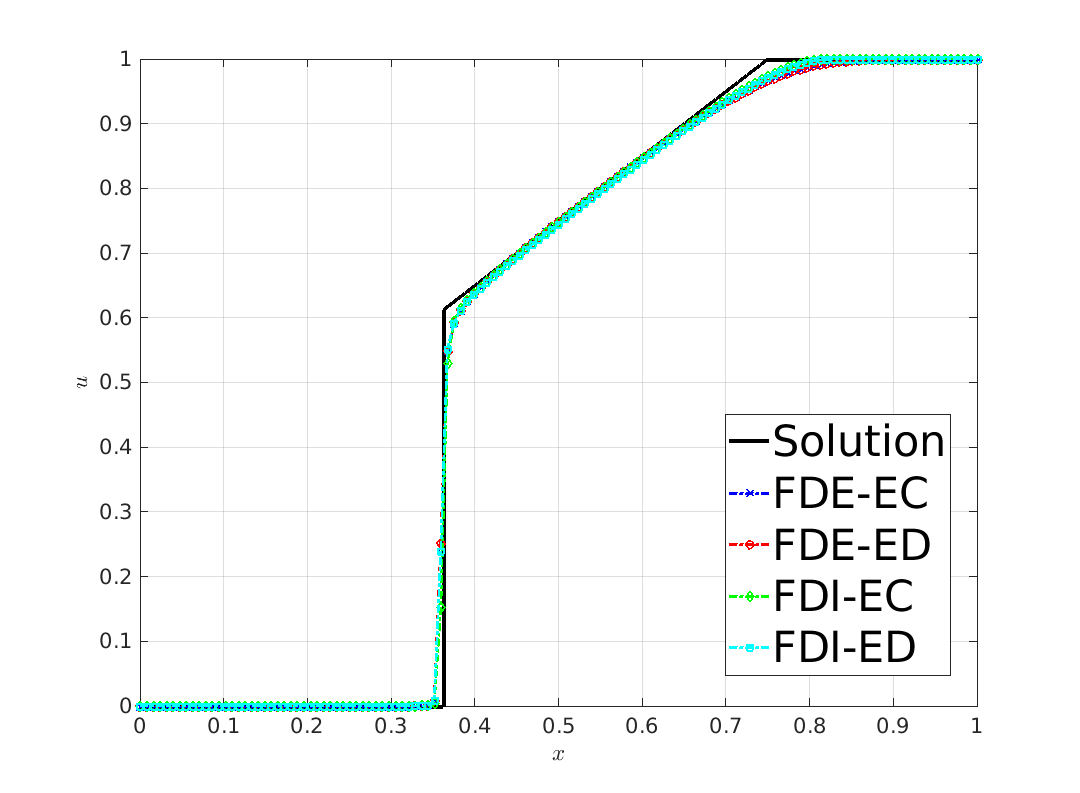}
\end{subfigure}
\begin{subfigure}[b]{0.49\textwidth}
\caption{RP2 at $t=2.0$, no fully discrete entropy fixes}
\includegraphics[width=\textwidth]{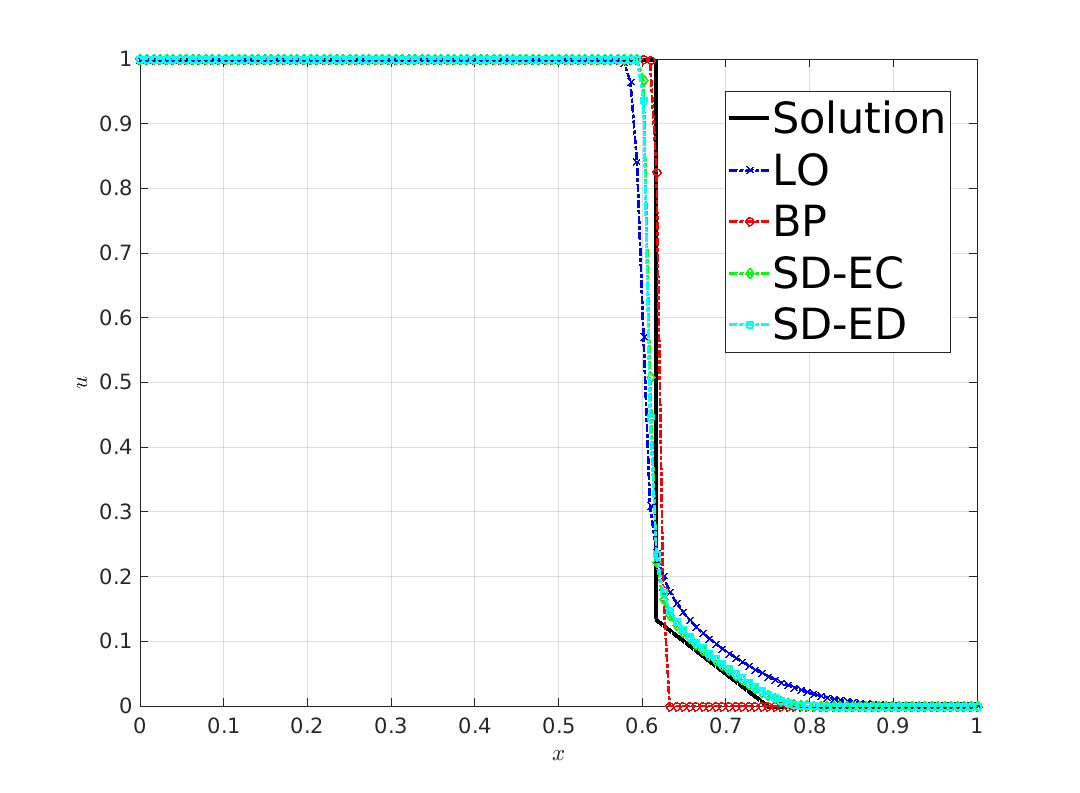}
\end{subfigure}
\begin{subfigure}[b]{0.49\textwidth}
\caption{RP2 at $t=2.0$ with fully discrete entropy fixes}
\includegraphics[width=\textwidth]{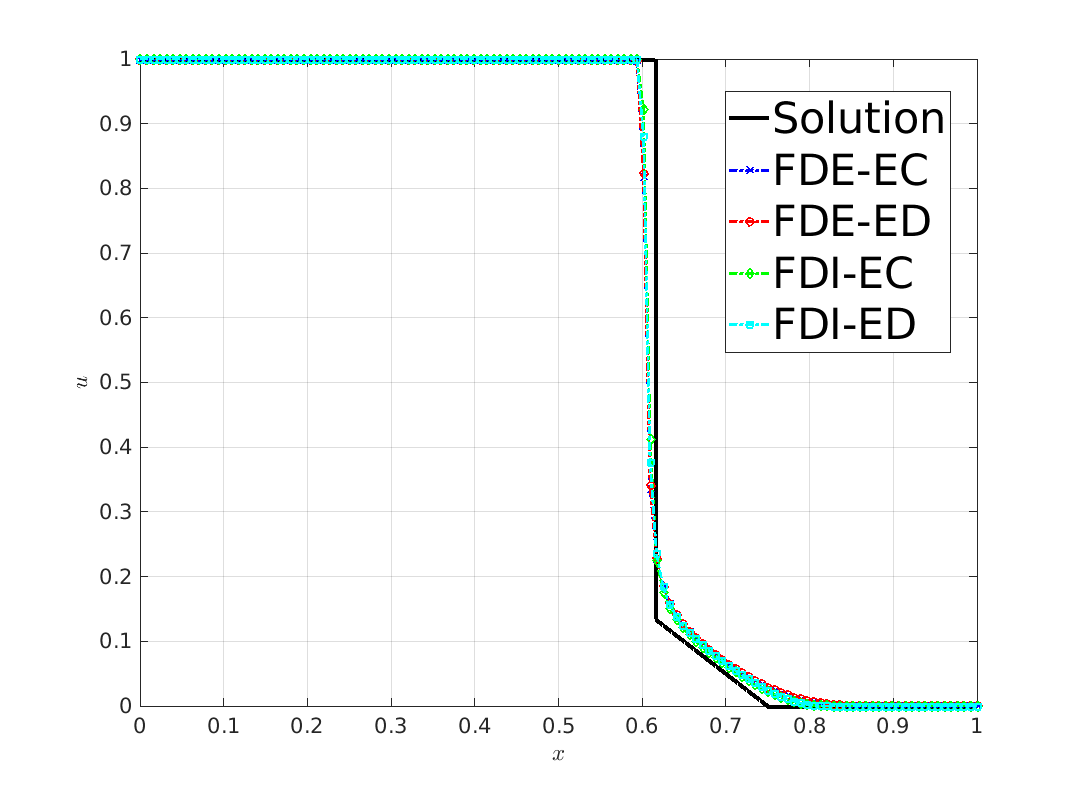}
\end{subfigure}
\caption{1D KPP problem,  numerical solutions of RP1 and RP2  calculated 
  with LO, BP, SD, and FD schemes on a uniform mesh of $\mathbb{P}_1$
 finite elements using $h=\frac{1}{128}$ and
 $\Delta t = 5\cdot 10^{-3}$. The acronyms used in the caption and
 legends are defined at the beginning of Section~\ref{sec:results}.}
\label{fig:kpp1D}
\end{figure}

We solve both Riemann problems numerically using the LO, BP, SD, and FD versions of the AFC scheme. We also vary the definition of the entropy production bound (EC vs. ED) and the type of the fully discrete fix (FDE vs. FDI) to study how these choices affect the entropy stability properties of the methods under investigation.
All profiles shown in \cref{fig:kpp1D} were computed on a uniform mesh with 128 cells using the fixed time step $\Delta t = 5 \cdot 10^{-3}$. In the RP1 and RP2 test alike, the BP scheme without entropy correction produces a wrong approximation in the post-shock region. The spurious plateaus in the red curves remain present as the mesh size and time step are refined. Hence, failure to perform an entropy fix inhibits convergence to vanishing viscosity solutions of RP1 and RP2.
Activation of the semi-discrete entropy fix is sufficient to cure this unsatisfactory behavior, while additional fully discrete fixes seem to be unnecessary in this example.
Another interesting observation is that flux-corrected solutions of both Riemann problems are rather insensitive to the choice of the entropy bounds, as the EC and ED curves are almost indistinguishable.

\subsection{Two-dimensional KPP equation}

The two-dimensional KPP problem \cite{guermond2016,guermond2017,kurganov2007} is a particularly challenging test for high-order numerical schemes. Equation \eqref{ibvp-pde} with the nonconvex flux function
\beq
\mathbf{f}(u)=(\sin(u),\cos(u))
\eeq
is solved in the computational domain
$\Omega_h=(-2,2)\times(-2.5,1.5)$ using the initial
condition
\beq \label{kpp_init}
u_0(x,y)=\begin{cases}
\frac{7\pi}{2} & \mbox{if}\quad \sqrt{x^2+y^2}\le 1,\\
\frac{\pi}{4} & \mbox{otherwise}.
\end{cases}
\eeq
The entropy flux corresponding to the square entropy
$\eta(u)=\frac{u^2}2$ is given by
$$\mathbf{q}(u)=(u\sin(u)+\cos(u),u\cos(u)-\sin(u)).$$
The maximum wave speed is bounded by
$\lambda=1$. We use this value in formula
\eqref{gms} for the artificial viscosity coefficients $d_{ij}$.
More accurate estimates can be found in \cite{guermond2017}.

\begin{figure}[ht!]
  
\begin{minipage}[t]{0.33\textwidth}
\centering (a) HO

\includegraphics[width=\textwidth,trim=0 0 0 20,clip]{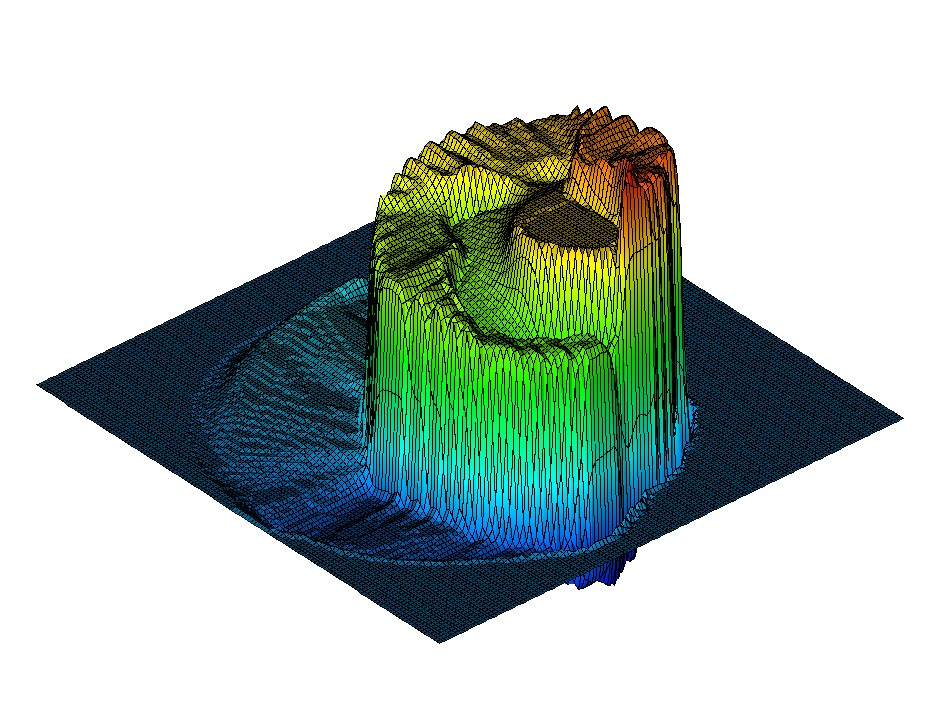}

\end{minipage}%
\begin{minipage}[t]{0.33\textwidth}

\centering (b) LO

\includegraphics[width=\textwidth,trim=0 0 0 50,clip]{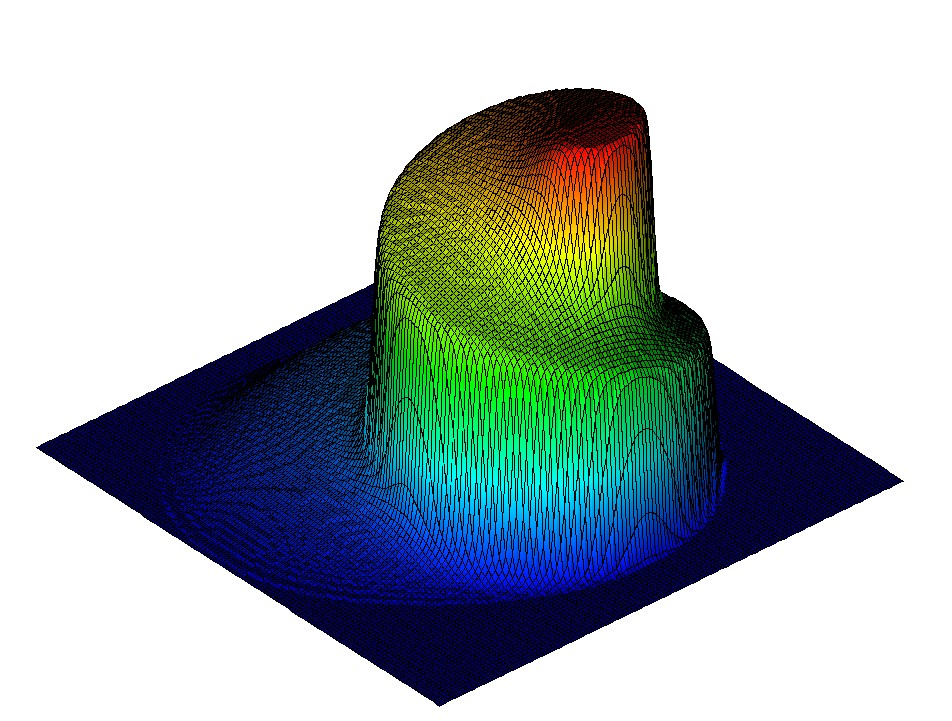}

\end{minipage}%
\begin{minipage}[t]{0.33\textwidth}

\centering (c) BP

\includegraphics[width=\textwidth,trim=0 0 0 20,clip]{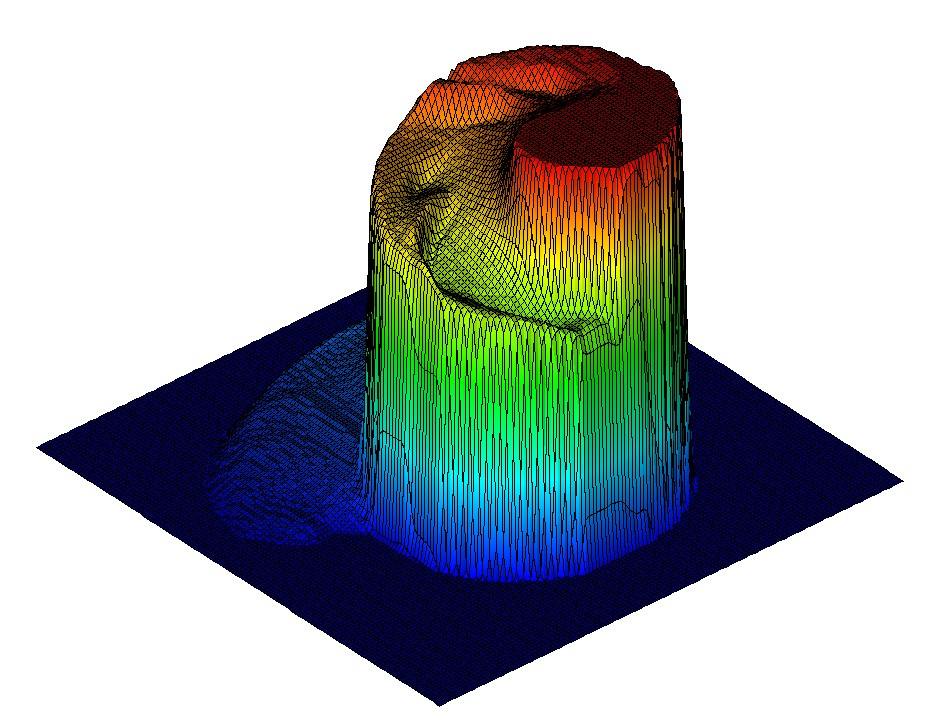}

\end{minipage}

\vskip0.25cm

\caption{KPP problem with initial condition \eqref{kpp_init}, (a) high-order, (b) low-order, and
 (c) flux-corrected  $\mathbb{Q}_1$ approximations
  at $t=1.0$ calculated on a uniform mesh using
  $h=\frac{1}{128}$, $\Delta t=10^{-3}$ without entropy fixes.}

\label{fig:KPP2D-1}
\vskip0.5cm

\begin{minipage}[t]{0.25\textwidth}
\centering (a) SD-EC

\includegraphics[width=\textwidth,trim=0 0 0 20,clip]{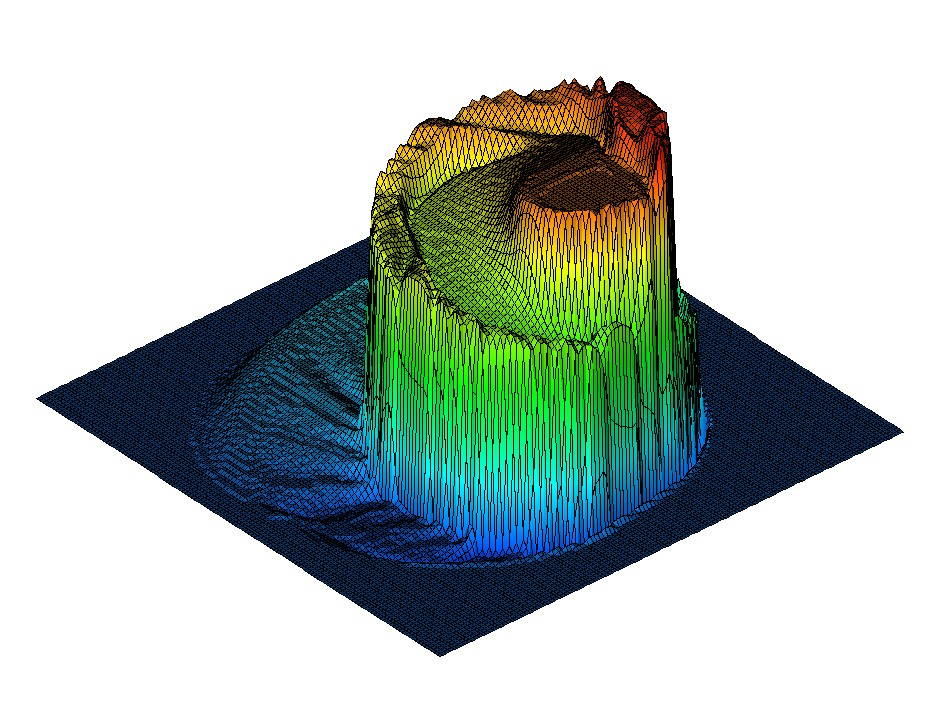}

\end{minipage}%
\begin{minipage}[t]{0.25\textwidth}

\centering (b) FDE-EC

\includegraphics[width=\textwidth,trim=0 0 0 50,clip]{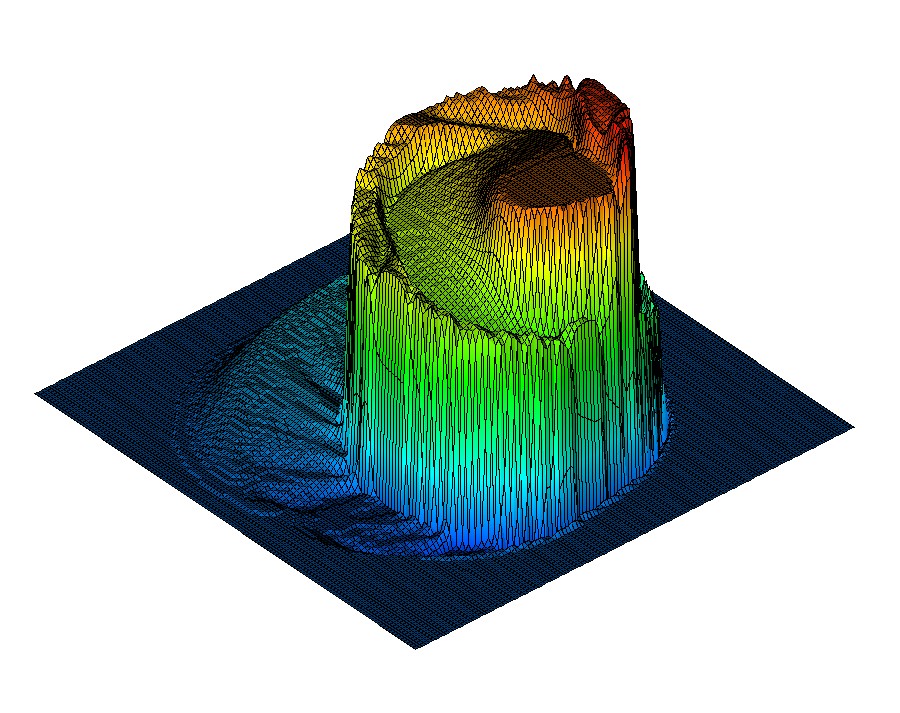}

\end{minipage}%
\begin{minipage}[t]{0.25\textwidth}
\centering (c)  SD-ED

\includegraphics[width=\textwidth,trim=0 0 0 20,clip]{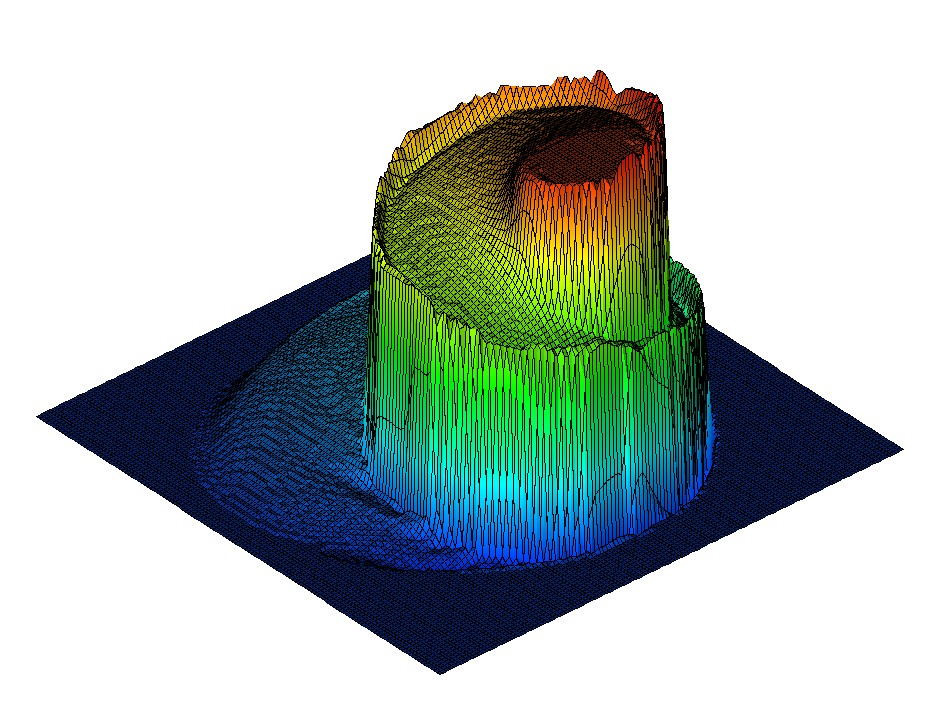}

\end{minipage}%
\begin{minipage}[t]{0.25\textwidth}

\centering (d) FDE-ED

\includegraphics[width=\textwidth,trim=0 0 0 50,clip]{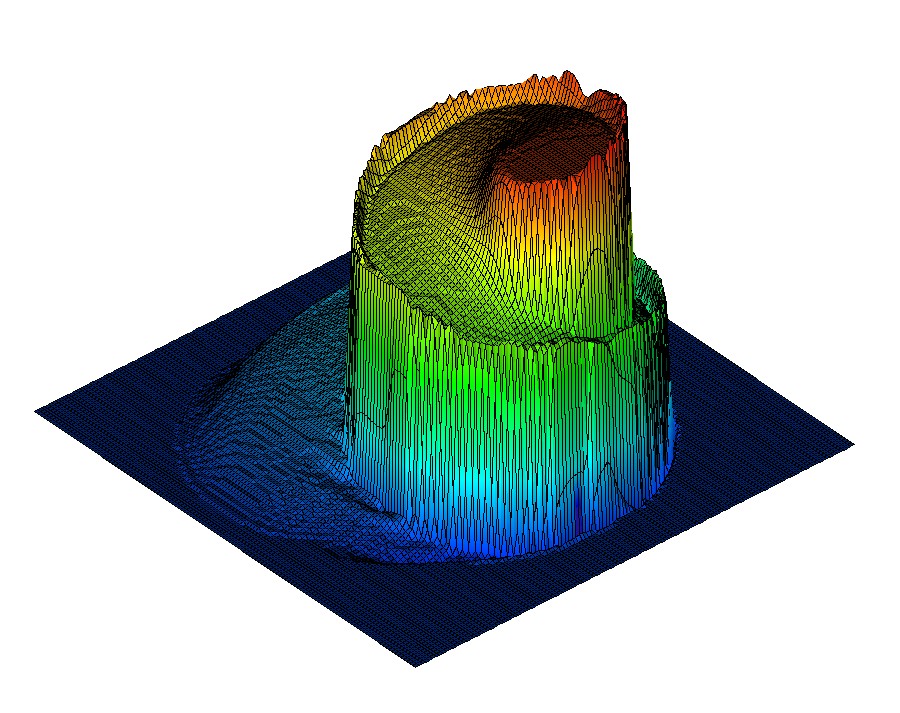}

\end{minipage}

\vskip0.25cm

\caption{KPP problem with initial condition \eqref{kpp_init}, flux-corrected $\mathbb{Q}_1$ approximations
  at $t=1.0$ calculated on a uniform mesh using
  $h=\frac{1}{128}$, $\Delta t=10^{-3}$
  and entropy stabilization without BP limiting.}

\label{fig:KPP2D-2}

\end{figure}

\begin{figure}[ht!]
  
\begin{minipage}[t]{0.33\textwidth}
\centering (a) SD-EC

\includegraphics[width=\textwidth,trim=0 0 0 20,clip]{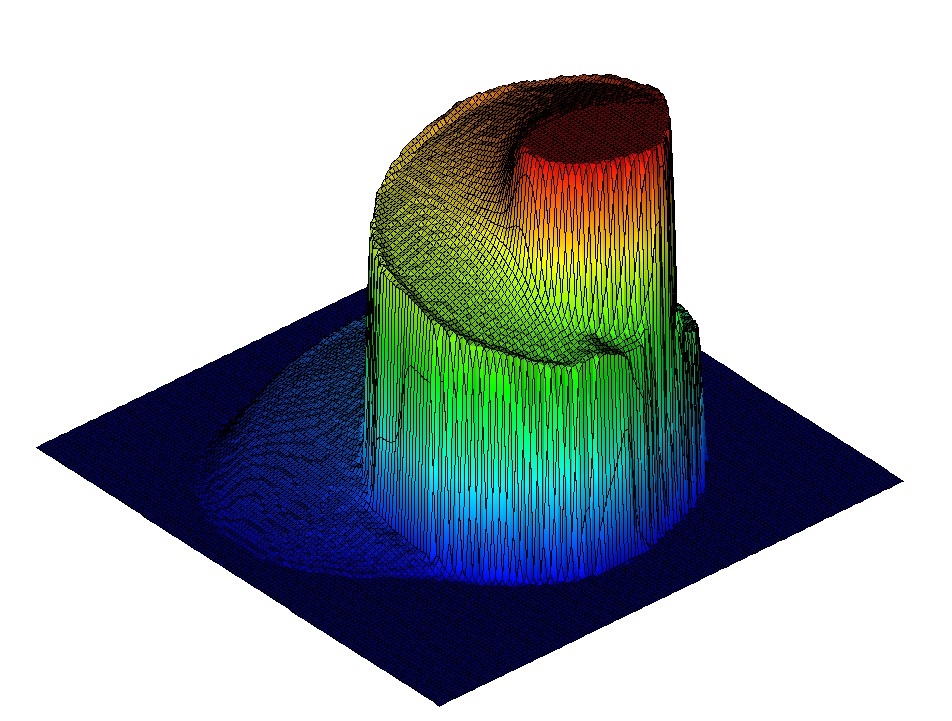}

\end{minipage}%
\begin{minipage}[t]{0.33\textwidth}

\centering (b) FDE-EC

\includegraphics[width=\textwidth,trim=0 0 0 50,clip]{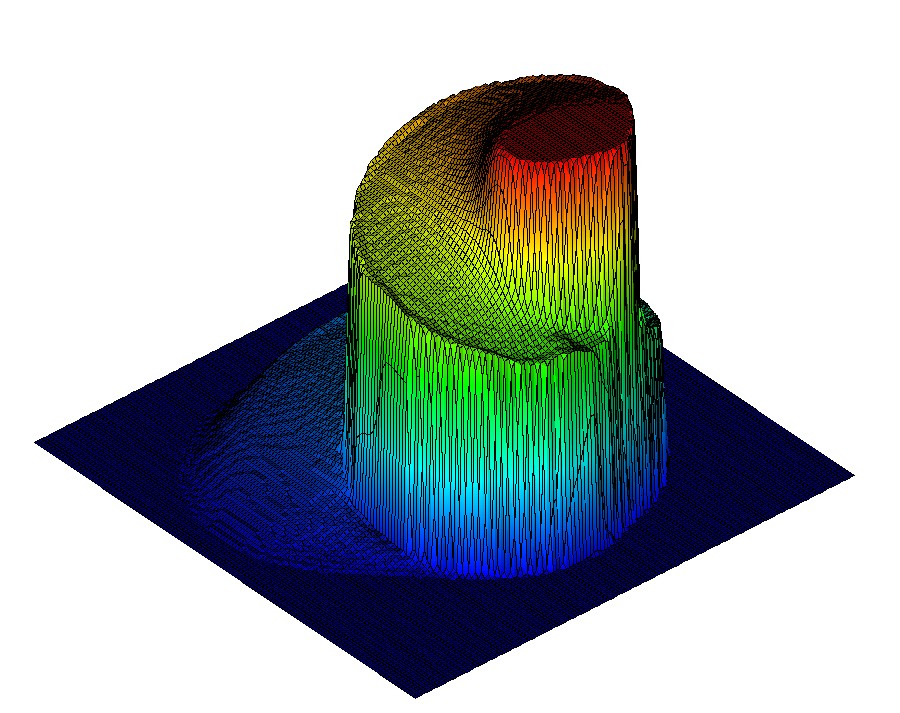}

\end{minipage}%
\begin{minipage}[t]{0.33\textwidth}

\centering (c) FDI-EC

\includegraphics[width=\textwidth,trim=0 0 0 20,clip]{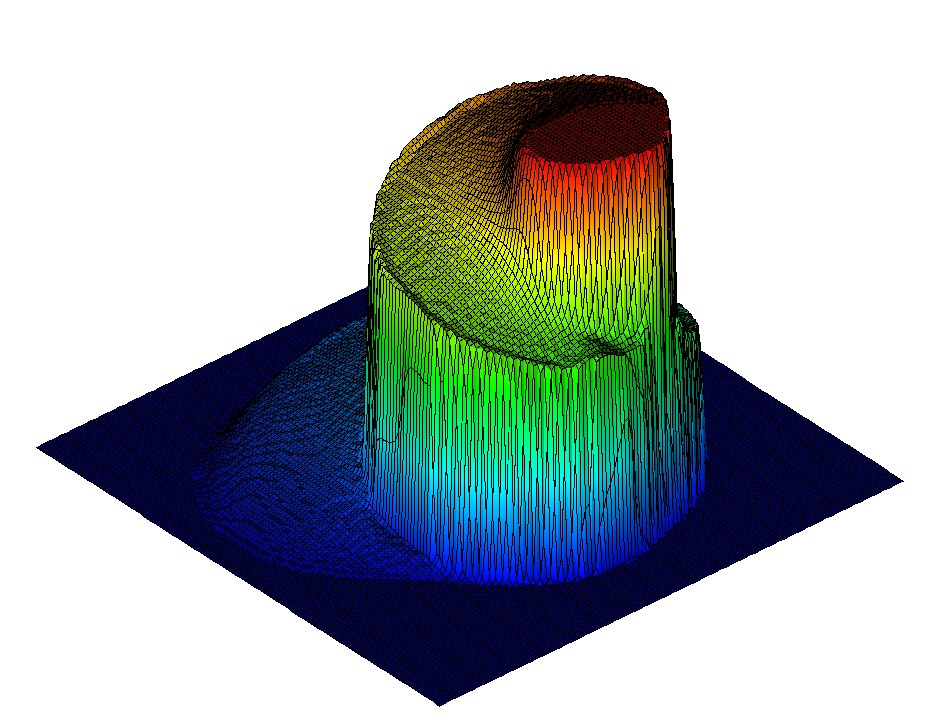}

\end{minipage}

\vskip0.25cm

\begin{minipage}[t]{0.33\textwidth}
\centering (d) SD-ED

\includegraphics[width=\textwidth,trim=0 0 0 20,clip]{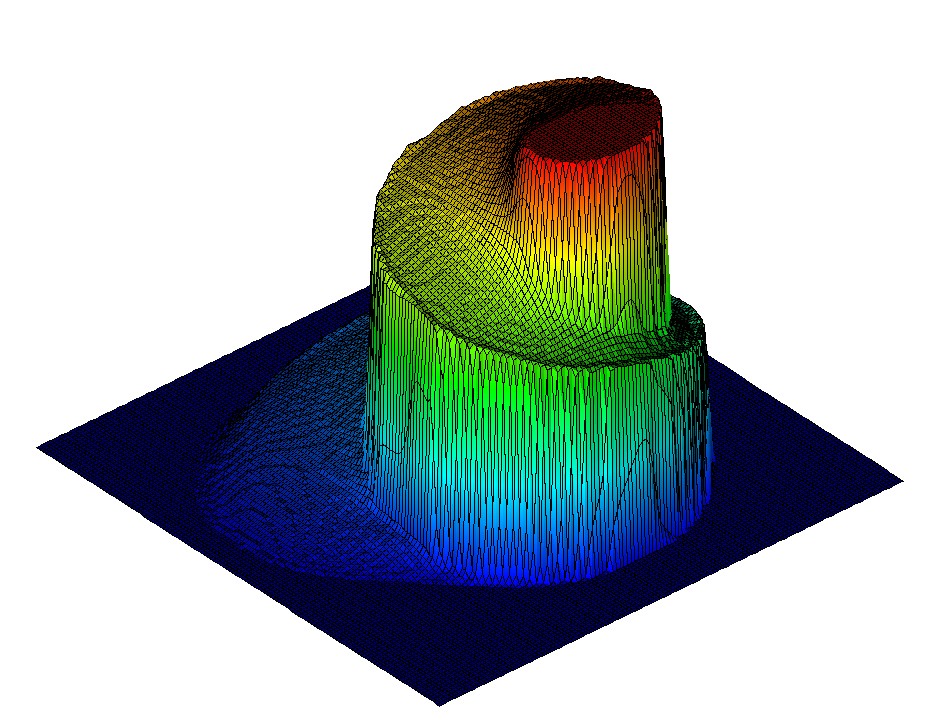}

\end{minipage}%
\begin{minipage}[t]{0.33\textwidth}

\centering (e) FDE-ED

\includegraphics[width=\textwidth,trim=0 0 0 50,clip]{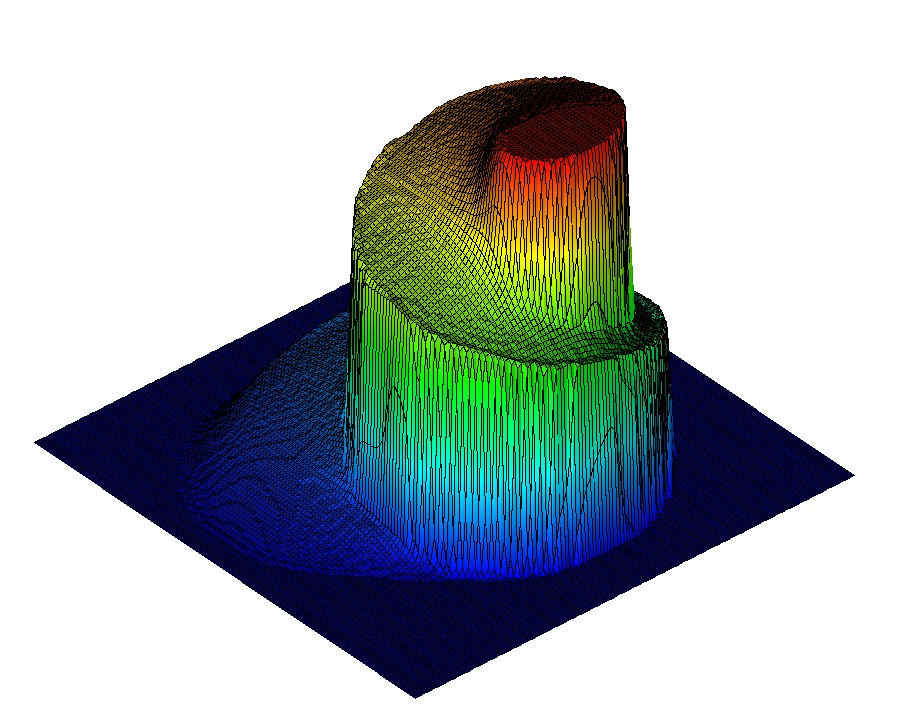}

\end{minipage}%
\begin{minipage}[t]{0.33\textwidth}

\centering (f) FDI-ED

\includegraphics[width=\textwidth,trim=0 0 0 20,clip]{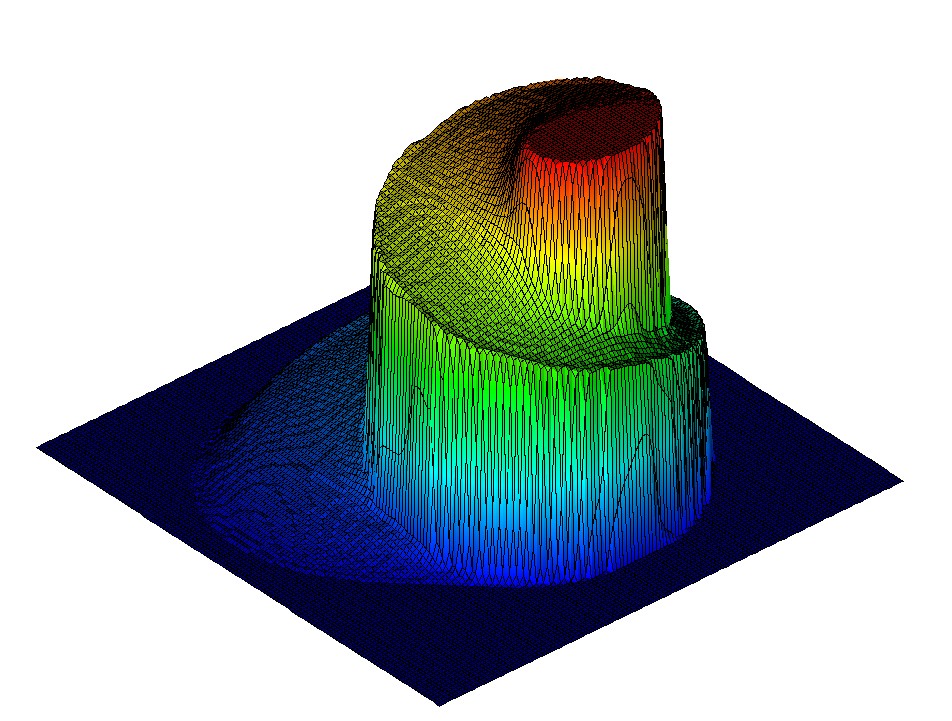}

\end{minipage}
\vskip0.25cm

\caption{KPP problem with initial condition \eqref{kpp_init}, flux-corrected $\mathbb{Q}_1$ approximations
  at $t=1.0$ calculated on a uniform mesh using
  $h=\frac{1}{128}$, $\Delta t=10^{-3}$
  and entropy stabilization with BP limiting.}

\label{fig:KPP2D-3}

\end{figure}

The entropy solution of 
the KPP problem exhibits a two-dimensional rotating wave structure.
The main challenge of this test is to avoid
convergence to wrong weak solutions. 
All results displayed in Figs~\ref{fig:KPP2D-1}--\ref{fig:KPP2D-3}
were calculated on a uniform mesh of $128\times 128$ bilinear
finite elements using the time step $10^{-3}$ and the final time
$t=1.0$. The oscillatory and qualitatively incorrect numerical
solution shown in Fig.~\ref{fig:KPP2D-1}(a) was produced by the
high-order baseline scheme, i.e., \eqref{semi-afc} with
$\alpha_{ij}=1$ and $f_{ij}^*=f_{ij}$ defined by \eqref{target}.
The low-order approximation shown in Fig.~\ref{fig:KPP2D-1}(b)
was obtained with $\alpha_{ij}=0$. It reproduces
the rotating wave structure of the entropy solution correctly
(cf. \cite{guermond2017,kuzmin2020g}) but is very diffusive.
The result of bound-preserving MCL flux correction without
entropy fixes is presented in Fig.~\ref{fig:KPP2D-1}(c).
There are no spurious oscillations but the AFC scheme converges to
a wrong weak solution.

In the next series of numerical experiments, we apply the entropy
correction factors $\alpha_{ij}$ to the unconstrained target fluxes
$f_{ij}^*=f_{ij}$, i.e., we deactivate the bound-preserving (BP) flux
limiter to study the effect of limiter-based entropy stabilization
separately. The results are shown in Fig.~\ref{fig:KPP2D-2}. It
can be seen that entropy-conservative (EC) fixes are insufficiently
dissipative, while their entropy-dissipative (ED) counterparts
perform well. The numerical solutions presented
in Fig. \ref{fig:KPP2D-3} demonstrate that BP flux limiting
eliminates undershoots and overshoots but convergence
to a wrong weak solution is possible if the EC bound
is used in the entropy fix. The
nonoscillatory approximations
produced by the BP-ED versions of the SD, FDE, and FDI
schemes preserve
the wave structure of the entropy-stable LO result
and are less diffusive. In this example, no significant differences
are observed between the outcomes of semi-discrete and fully
discrete fixes based on the same definition of $Q_{ij}$.

To study the convergence behavior of different algorithms for
a nonlinear problem with a smooth exact solution, we replace
the discontinuous initial condition \eqref{kpp_init} by
\beq \label{kpp_init_smooth}
u_0(x,y)=\begin{cases}
\frac{\pi}{4}\left(1+\frac{1}{20}\left(1+\cos(\pi\sqrt{x^2+y^2}\right)\right)
& \mbox{if}\quad \sqrt{x^2+y^2}\le 1,\\
\frac{\pi}{4} & \mbox{otherwise}
\end{cases}
\eeq
and determine the experimental order of convergence (EOC)
for $h=\frac{1}{256}$ using the formula
$$
\mathrm{EOC}=\frac{\log\left(\frac{\|u_{4h}-u_{2h}\|}{\|u_{2h}-u_h\|}\right)}{\log 2}.
$$

\begin{table}[b!]
\centering

\tabcolsep0.4cm
\renewcommand{\arraystretch}{1.1}

\begin{tabular}{c|cccccc}
\hline
   & LO & HO & BP & SD & FDE & FDI    \\\hline
  $\mathrm{EOC}_{L^1}$ & 0.75 & 2.28 & 2.39 & 2.40 & 1.28 & 2.13\\
  $\mathrm{EOC}_{L^2}$ & 0.71 & 2.06 & 2.25 & 2.30 & 1.24 & 2.11
\\\hline
\end{tabular}

\vskip0.25cm

\caption{KPP problem with initial condition \eqref{kpp_init_smooth},
  $L^1$ and $L^2$ convergence rates.}
\label{kpp_eoc}
\end{table}

In the process of mesh refinement, we keep the ratio
$\frac{\Delta t}{h}=0.256$ fixed.
The rates of convergence w.r.t. discrete $L^1$ and  $L^2$ norms
are reported in Table \ref{kpp_eoc}. In this experiment, all
entropy fixes use BP prelimiting and the ED bound defined by \eqref{Qij-ED}.
The EOC of the unconstrained HO scheme is slightly higher than
2 due to superconvergence of consistent-mass approximations
using linear or bilinear finite elements. The BP, SD, and FDI
algorithms preserve second-order accuracy. The first-order
convergence behavior of the FDE version is caused by the
failure of stagewise explicit fixes to exploit cancellation
of entropy production and dissipation in high-order Runge--Kutta schemes
(see Remark \ref{rmk:fde}). If the time step is chosen so that 
$\frac{\Delta t}{h^2}=128^2\cdot 10^{-3}$ remains fixed, the
 FDE version delivers $\mathrm{EOC}_{L^1}$=1.95 and $\mathrm{EOC}_{L^2}$=1.85.
 To achieve second-order accuracy 
without using $\Delta t=\mathcal O(h^2)$, the entropy correction factors 
of the FDE fix could be redefined as $$\alpha_{ij}=\max\{1-
\max\{\beta_i,\beta_j\},\min\{R_i,R_j\}\},$$
where $\beta_i\in[0,1]$ is a smoothness indicator such as
the entropy residual sensor used in \cite{kuzmin2020f} to
stabilize the target flux as originally proposed in
\cite{guermond2018}. The FDI fix can also be configured
to use $\alpha_{ij}\ge 1-\max\{\beta_i,\beta_j\}$. In this
way, the use of fully discrete fixes can be restricted
to subdomains with large entropy residuals and the cost
of calculating $\alpha_{ij}$ can be reduced.

To compare our new approaches with an existing FD fix, we used
the method of Berthon et al.\ \cite{berthon2021,berthon2020} to construct
localized artificial viscosities $d_{ij}^{\mathrm{FD}}$ for the
entropy-corrected fluxes $\alpha_{ij}f_{ij}^*+d_{ij}^{\mathrm{FD}}
(u_j-u_i)$, where $\alpha_{ij}$ is the correction factor of the
semi-discrete entropy fix \eqref{alpha}. In our AFC notation,
the formula
for $d_{ij}^{\mathrm{FD}}$ is given by (cf. \cite[eq.(12)--(14)]{berthon2020})
 \beq\label{barfix}
 d_{ij}^{\mathrm{FD}}=\begin{cases}
 -\frac{P_{ij}}{2D_{ij}} & \mbox{if}\ P_{ij}D_{ij}<0,\\
 0 & \mbox{otherwise},
 \end{cases}
 \eeq
\begin{align*}
 P_{ij}&=d_{ij}[\eta(u_j^*)+\eta(u_i^*)-\eta(u_j)+\eta(u_i)]
 +(\mathbf{q}_j-\mathbf{q}_i)\cdot\mathbf{c}_{ij},\\
 D_{ij}&=2\eta\left(\frac{u_j+u_i}2\right)-\eta(u_j)-\eta(u_i),\\ 
u_j^*&=u_j-\frac{(\mathbf{f}_j-\mathbf{f}_i)\cdot\mathbf{c}_{ij}
+d_{ij}(u_j-u_i)+\alpha_{ij}f_{ij}^*}{2d_{ij}},\\
u_i^*&=u_i-\frac{(\mathbf{f}_j-\mathbf{f}_i)\cdot\mathbf{c}_{ij}
-d_{ij}(u_j-u_i)-\alpha_{ij}f_{ij}^*}{2d_{ij}}.
\end{align*}
Note that $u_j^*+u_i^*=2\bar u_{ij}$ and $D_{ij}\le 0$ since $\eta(u)$
is convex. The analysis in \cite{berthon2021,berthon2020} is restricted
to 1D and assumes that the same global value of $d_{ij}^{\mathrm{FD}}$
is used
for all fluxes. Hence, the multidimensional
generalization \eqref{barfix} may fail to enforce fully discrete entropy
stability. Nevertheless, it provides an interesting
alternative to the limiter-based FDE fix.

Lemma 2 in \cite{berthon2020} ensures that $d_{ij}^{\mathrm{FD}}=\gamma d_{ij}$,
where $\gamma=\mathcal O(1)$, for target fluxes of the form
$f_{ij}=d_{ij}(u_i-u_j)$. However, the ones defined by \eqref{target}
do not necessarily vanish if $u_i=u_j$. As a consequence, $\gamma$ can
become unbounded or large enough to violate the CFL condition.
For that reason, we use $d_{ij}$ defined by \eqref{gms} as upper
bound for $d_{ij}^{\mathrm{FD}}$ in the consistent-mass version.
The results for the KPP problem with the
discontinuous initial condition \eqref{kpp_init} are similar to those
obtained with the FDE fix (not shown here). As expected, the
artificial viscosity method also exhibits first-order
convergence ($\mathrm{EOC}_{L^1}=0.98$, $\mathrm{EOC}_{L^2}$=0.95) in the smooth
KPP test with the initial condition \eqref{kpp_init_smooth}
and $\frac{\Delta t}{h}=0.256$ refinements.

\subsection{One-dimensional shallow water equations}\label{sec:dam}

Having studied the KPP problem in the one- and two-dimensional setting, we now move on to the case of 1D systems of conservation laws.
First, we consider the shallow water equations. In the case of a flat bottom topography, this hyperbolic system reads
\begin{align*}
\frac{\partial}{\partial t} \begin{bmatrix}
h\\
hv 
\end{bmatrix}  + \frac{\partial}{\partial x}
\begin{bmatrix}
hv \\
hv^2 + \frac g 2 h^2
\end{bmatrix} = 0.
\end{align*}
Here $h$ is the total water height, $v$ is the depth-integrated velocity, and $g$ is the gravitational acceleration, which we set equal to 1.
An entropy-entropy flux pair and the corresponding entropy potential for the shallow water system are given by \cite{chen2017,fjordholm2011}
\begin{align*}
\eta(u) = \frac 1 2 (hv^2 + g h^2),\qquad q(u) = \frac 1 2 h v^3 + gh^2 v, \qquad \psi(u) = \frac 1 2 gh^2 v.
\end{align*}

We test the numerical behavior of the proposed flux correction methods for a wet dam break example corresponding to the initial condition
\begin{align}\label{eq:swe-init}
u_0(x) = \begin{cases}
(1,0) & \mbox{if } x < 0, \\ (0.1,0) & \mbox{otherwise.}
\end{cases}
\end{align}
In this Riemann problem, two water columns of different height are initially separated by a dam, which is removed at the start of the simulation.
After the dam break, the higher water column expands into a rarefaction fan, and a shock front propagates into the lower water level region. The height of the plateau between the shock and the rarefaction wave is given by $c_m^2/g$, where $c_m$ is a root of a sixth-degree polynomial.
For the particular initial condition \eqref{eq:swe-init}, we use the value $c_m \approx 0.6294$ to approximate the analytical solution of this problem as in \cite{delestre-arxiv} (the expressions in the journal version of this reference are incorrect).

We simulate the dam break in the computational domain $\Omega = (-0.5,0.5)$ and impose reflecting wall boundary conditions at $x=\pm 0.5$. The final time is $t=0.3$. To compare the accuracy of individual approaches, we run simulations on a sequence of successively refined uniform meshes using the time step $\Delta t = 0.25h$. The scalar quantity of interest
\begin{align*}
e_1(t) = \|u(t)-u_h(t)\|_{L^1(\Omega)^2}
\end{align*}
is defined as the sum of $L^1(\Omega)$ errors in the conserved variables.
It measures the accuracy of a given approximation at $t\ge 0$.
The values of $e_1(0.3)$ and the corresponding EOCs are summarized in \cref{tab:swe-gt,tab:swe-rt}.  Because of the discontinuity in \eqref{eq:swe-init}, first-order accuracy of the BP scheme equipped with the Galerkin target flux is optimal.  The semi-discrete entropy fix using either \eqref{Qij-EC} or \eqref{Qij-ED} does not degrade the rate of convergence. The high-resolution GT versions of the BP, SD, and FD approaches are clearly superior to the low-order method, which is significantly more diffusive.
All numerical results obtained with the Roe target flux are also less accurate than the corresponding Galerkin approximations.

\begin{figure}[ht!]
\centering
\begin{subfigure}[b]{0.49\textwidth}
\caption{GT-BP}
\includegraphics[width=\textwidth]{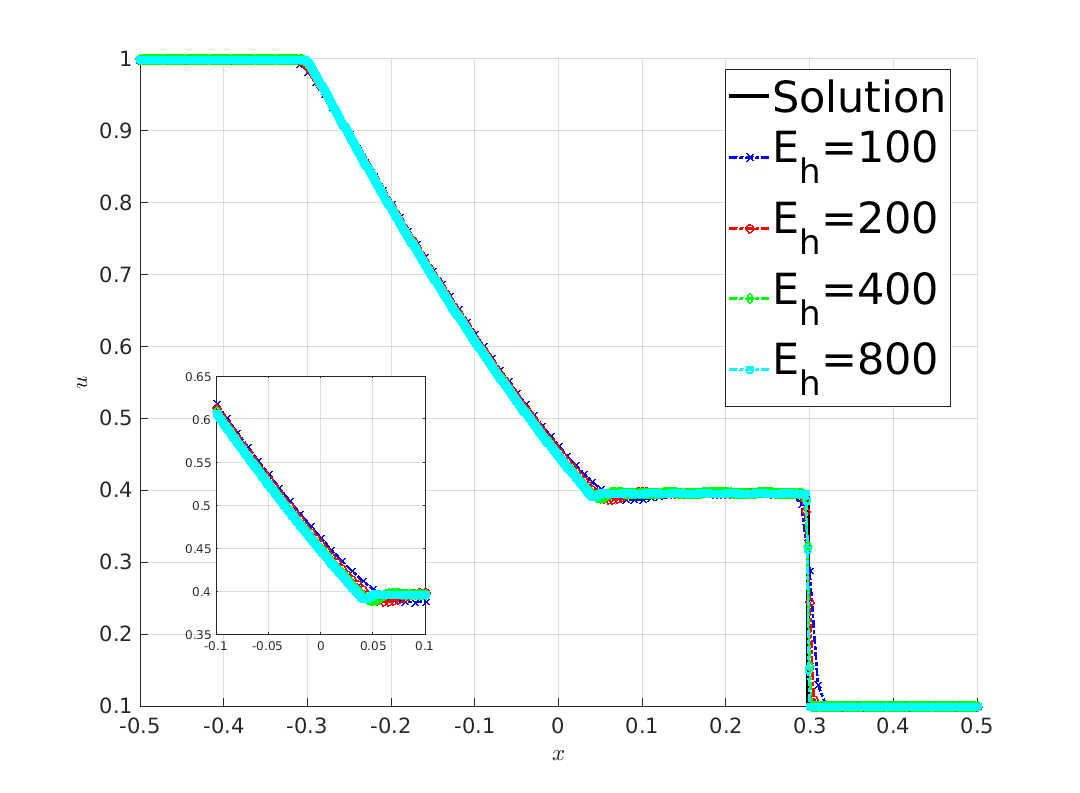}
\end{subfigure}
\begin{subfigure}[b]{0.49\textwidth}
\caption{RT-BP}
\includegraphics[width=\textwidth]{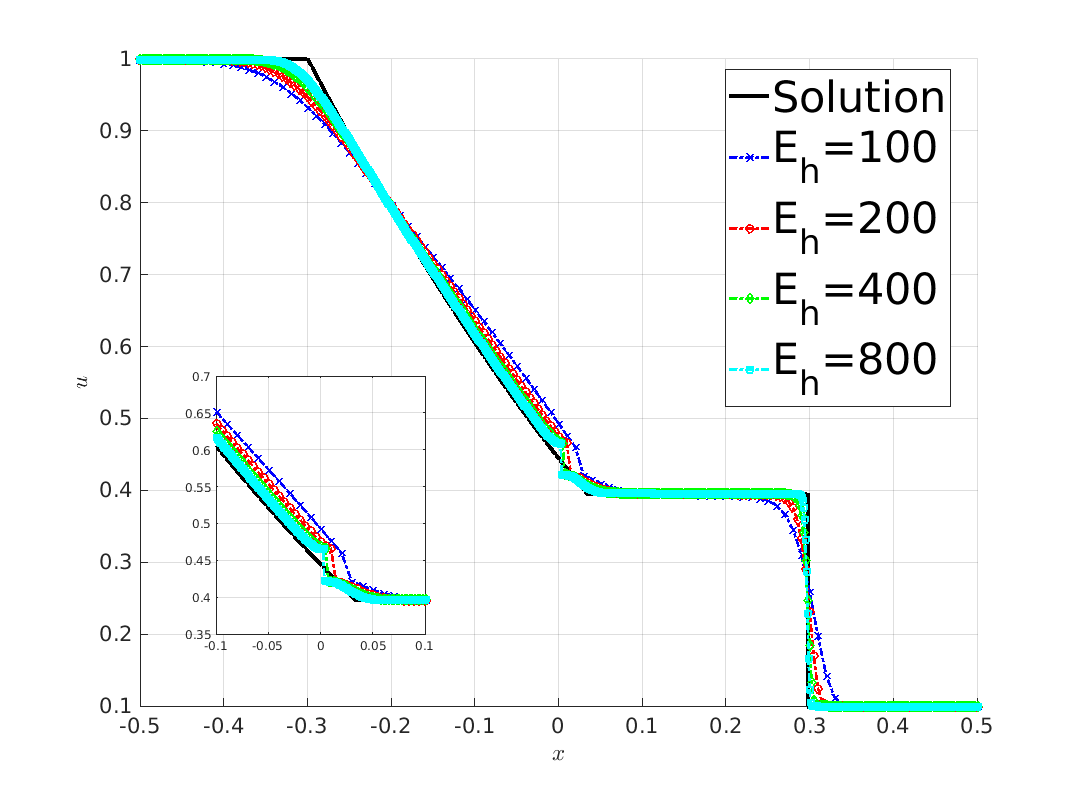}
\end{subfigure}
\begin{subfigure}[b]{0.49\textwidth}
\caption{RT-SD}
\includegraphics[width=\textwidth]{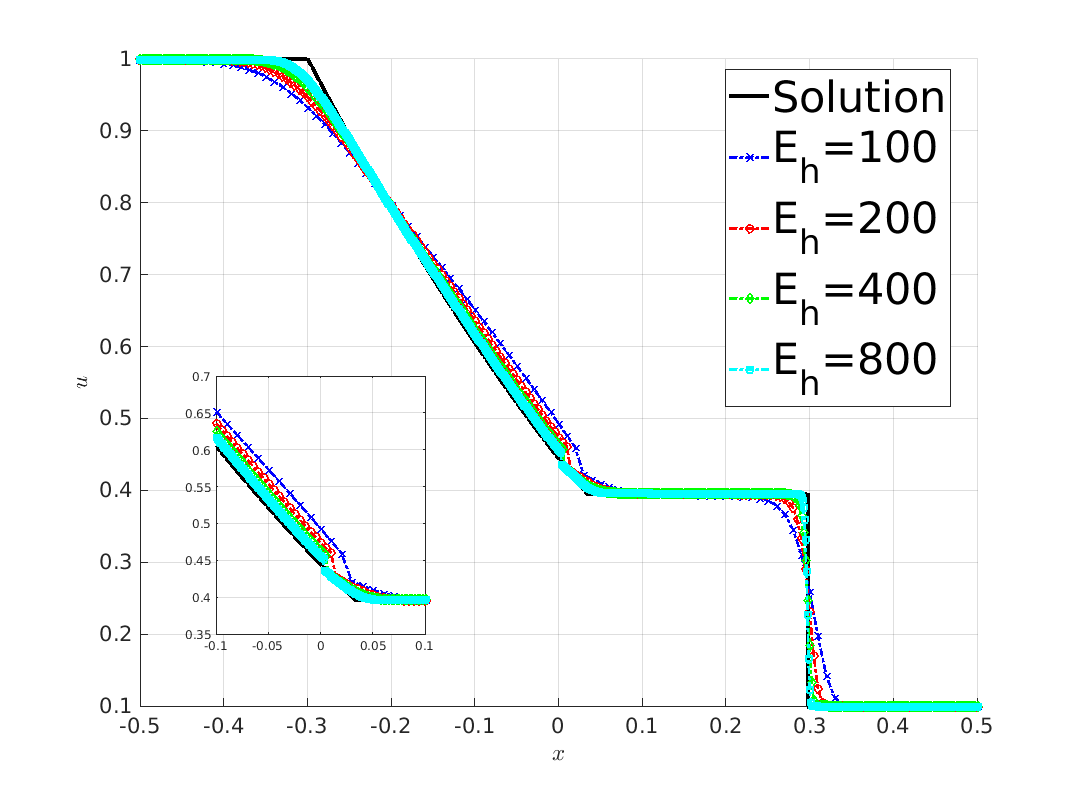}
\end{subfigure}
\begin{subfigure}[b]{0.49\textwidth}
\caption{RT-FDE}
\includegraphics[width=\textwidth]{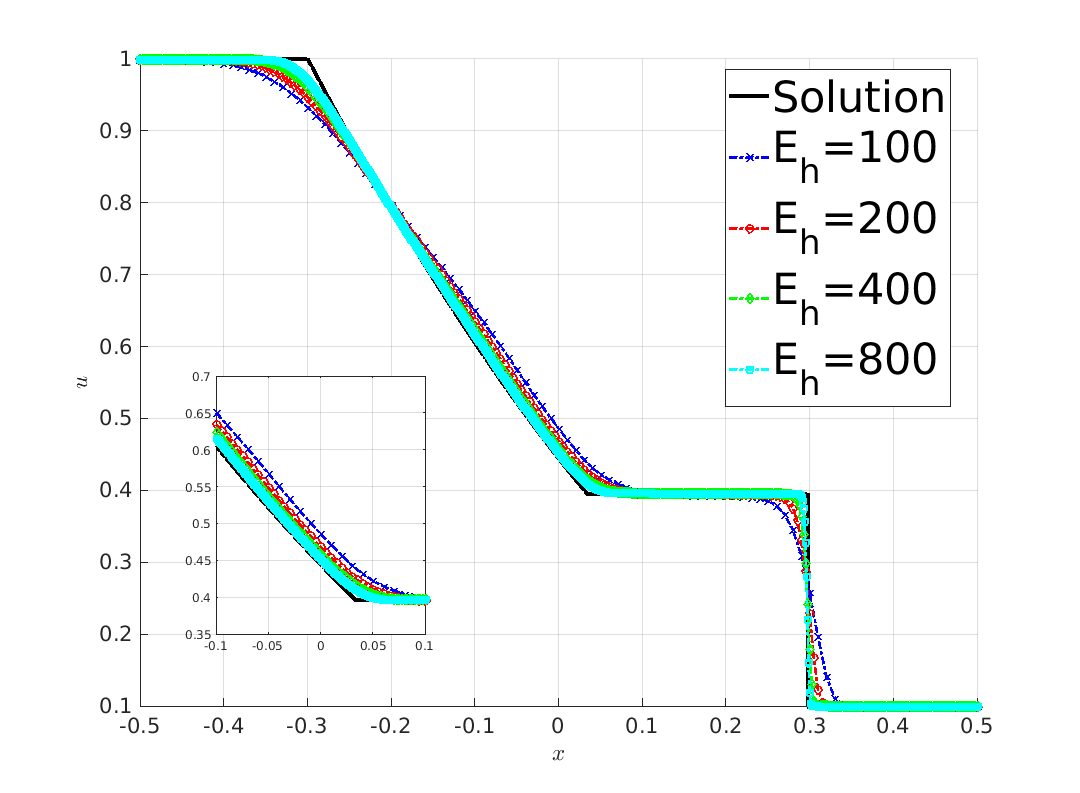}
\end{subfigure}
\caption{Vanishing viscosity solution and approximations to the water height of the wet dam break with initial condition \eqref{eq:swe-init}. Results at $t=0.3$ on four refinement levels using $\Delta t = 0.25h$.}\label{fig:dam}
\end{figure}

It is well known that Roe's scheme applied to the Euler equations can produce nonphysical stationary shocks at sonic points \cite{banks2009,chandrashekar2013,einfeldt1988,moschetta2000,toro2009}.
The same is true for the shallow water equations \cite{kemm2014,ketcheson-arxiv}.
In fact, an entropy-violating shock can form even in our simple example.
To show this and study the convergence behavior of GT-BP, RT-BP, RT-SD, and RT-FDE approximations, we performed a sequence of simulations with increased spatial and temporal resolution. The results for the water height are shown in \cref{fig:dam}.
Since the errors and EOCs in \cref{tab:swe-gt,tab:swe-rt} indicate that EC and ED versions of entropy fixes perform similarly for this benchmark, we only show the results obtained with entropy-conservative bounds.
The use of Galerkin target fluxes produces profiles that are in very good agreement with the analytical solution, while solutions for Roe target fluxes exhibit small nonphysical shocks within the rarefaction wave region.
The semi-discrete entropy fix reduces the magnitude of spurious jumps but only the additional fully discrete fix produces qualitatively correct profiles.

\begin{table}[ht!]
\centering
\begin{tabular}{c|cc|cc|cc|cc}
$E_h$ & LO & EOC & GT-BP & EOC & GT-SD-EC & EOC & GT-SD-ED & EOC \\
\hline
32  & 1.38E-01 &      & 5.99E-02 &      & 6.50E-02 &      & 6.57E-02 \\
64  & 8.43E-02 & 0.71 & 3.16E-02 & 0.92 & 3.42E-02 & 0.93 & 3.46E-02 & 0.92 \\
128 & 4.98E-02 & 0.76 & 1.61E-02 & 0.98 & 1.75E-02 & 0.97 & 1.77E-02 & 0.97 \\
256 & 2.91E-02 & 0.78 & 8.19E-03 & 0.97 & 8.88E-03 & 0.98 & 8.99E-03 & 0.98 \\
\end{tabular}
\caption{Wet dam break with initial condition \eqref{eq:swe-init}, convergence history for $e_1(0.3)$ and corresponding EOCs of low- and high-order schemes using Galerkin target fluxes.}\label{tab:swe-gt}
\end{table}
\begin{table}[ht!]
\centering
\begin{tabular}{c|cc|cc|cc|cc}
$E_h$ & RT-BP & EOC & RT-SD & EOC & RT-FD-EC & EOC & RT-FD-ED & EOC \\
\hline
32  & 1.33E-01 &      & 1.33E-01 &      & 1.33E-01 &      & 1.33E-01 \\
64  & 7.90E-02 & 0.75 & 7.90E-02 & 0.75 & 7.97E-02 & 0.74 & 7.97E-02 & 0.74 \\
128 & 4.57E-02 & 0.79 & 4.57E-02 & 0.79 & 4.62E-02 & 0.79 & 4.62E-02 & 0.79 \\
256 & 2.62E-02 & 0.80 & 2.61E-02 & 0.81 & 2.65E-02 & 0.80 & 2.65E-02 & 0.80 \\
\end{tabular}
\caption{Wet dam break with initial condition \eqref{eq:swe-init}, convergence history for $e_1(0.3)$ and corresponding EOCs high-order schemes using Roe target fluxes.}\label{tab:swe-rt}
\end{table}

\subsection{One-dimensional Euler equations of gas dynamics}

Another hyperbolic model of particular importance are the Euler equations of gas dynamics. In the 1D case, this system of conservation laws reads
\begin{align*}
\frac{\partial}{\partial t}\begin{bmatrix}
\rho\\
\rho v\\
\rho E
\end{bmatrix} + \frac{\partial}{\partial x}
\begin{bmatrix}
\rho v \\ \rho v^2 + p \\
(\rho E + p)v
\end{bmatrix} = 0,
\end{align*}
where $\rho$ is the fluid density, $v$ is the velocity, and $E$ is the specific total energy.
The pressure $p$ of a polytropic ideal gas is given by the equation of state
\begin{align*}
p(u) = (\gamma-1) \left( \rho E - \frac{\rho v^2}{2} \right)
\end{align*}
in which $\gamma$ denotes the adiabatic constant. For diatomic gases and
 air, $\gamma$ is approximately equal to $1.4$. We use this value in our
simulations. An entropy-entropy flux pair and the corresponding entropy potential for the Euler system are given by \cite{chen2017,harten1983}
\begin{align*}
\eta(u) =\frac{\rho s(u)}{1-\gamma},\qquad q(u) = \frac{\rho v s(u)}{1-\gamma}, \qquad \psi(u) = \rho v,
\end{align*}
where $s(u) = \log(p\rho^{-\gamma})$ is the specific entropy.
Note that $s(u)$ and $\eta(u)$ have opposite signs, since the physical entropy is concave while the mathematical one is convex.

\subsubsection{Shock tube problems}\label{sec:sod}

In the first numerical experiment for the one-dimensional Euler equations, we solve the classical
shock tube problem with Sod's \cite{sod1978} initial condition
\begin{align}\label{eq:sod}
u_0(x) = \begin{cases}
(1,0,2.5) & \mbox{if } x < 0.5, \\
(0.125,0, 0.25)& \mbox{otherwise}
\end{cases}
\end{align}
in the computational domain $\Omega=(0,1)$ with reflecting wall boundaries.
The analytical solution is qualitatively similar to that of the dam break
problem considered in \cref{sec:dam} but features a contact discontinuity
in addition to the shock and rarefaction.
To facilitate direct comparison with the results presented in \cite{kuzmin2020}, we stop simulations at $t=0.231$ and run them on a
uniform mesh with $h=\frac{1}{128}$ using the constant time 
step $\Delta t=10^{-3}$.

\begin{figure}[ht!]
\centering
\begin{subfigure}[b]{0.49\textwidth}
\caption{GT-SD-EC}
\includegraphics[width=\textwidth]{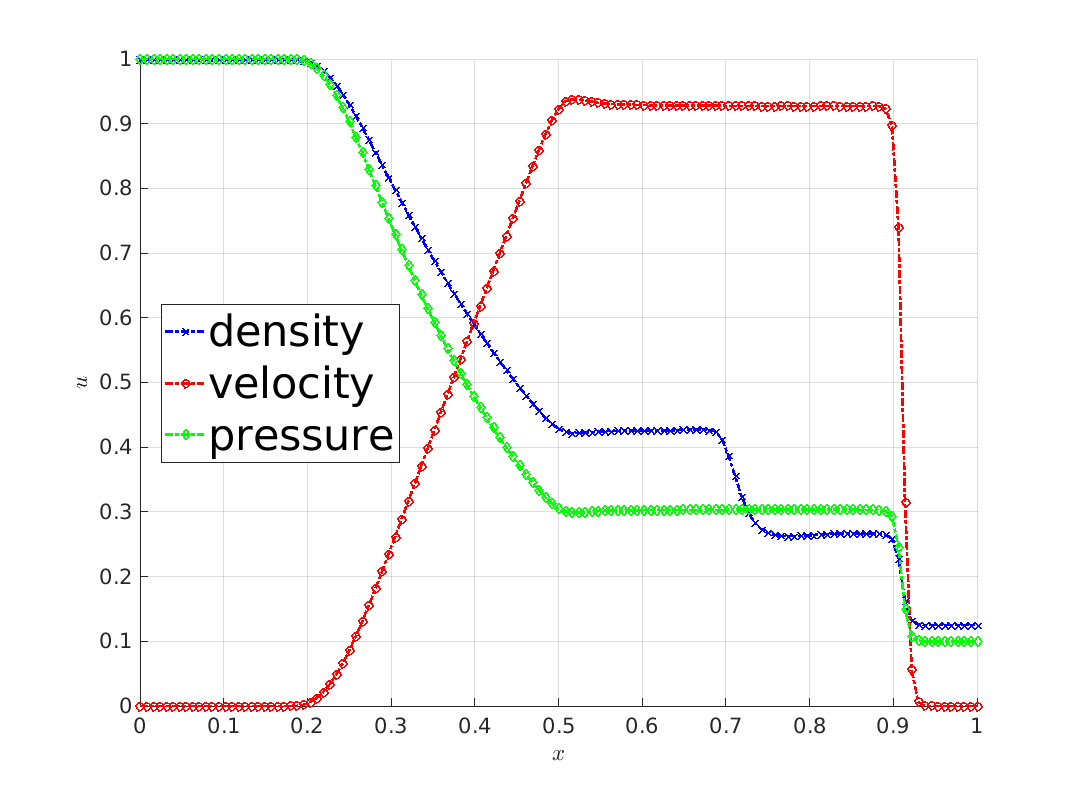}
\end{subfigure}
\begin{subfigure}[b]{0.49\textwidth}
\caption{RT-SD-EC}
\includegraphics[width=\textwidth]{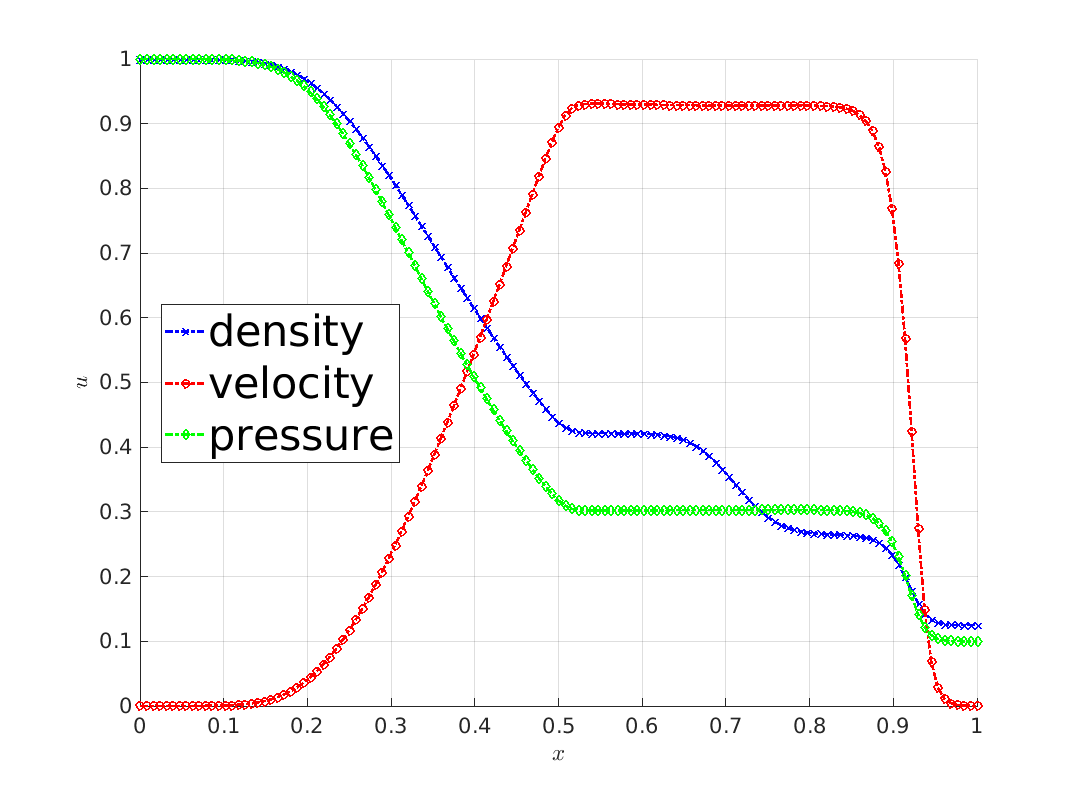}
\end{subfigure}
\caption{Sod's shock tube \cite{sod1978}, 
  primitive variables at $t=0.231$
  calculated on a uniform mesh with $h=\frac 1 {128}$ and $\Delta t= 10^{-3}$.
  The SD-EC scheme constrains the Galerkin or Roe target flux
  using bound-preserving convex limiting and
  the semi-discrete entropy fix with EC bounds
  \eqref{Qij-EC}.}\label{fig:sod}
\end{figure}

For this simple test problem, all high-resolution schemes based on the same
definition of the target flux produce similar results. That is why only
numerical solutions obtained with the GT and RT versions of the SD-EC
algorithm are compared in \cref{fig:sod}. The GT profiles are very similar
to the ones obtained without the semi-discrete entropy fix (cf. Fig.~4(a)
in \cite{kuzmin2020}) and less diffusive than the corresponding RT
approximations. In this example, all approaches produce qualitatively
correct results which converge to the entropy solution.

Let us now repeat the above shock tube experiment using the
initial condition
\begin{align}\label{eq:modsod}
u_0(x) = \begin{cases}
(1,0.75,89/32) & \mbox{if } x < 0.3, \\
(0.125,0, 0.25)& \mbox{otherwise}.
\end{cases}
\end{align}
The Riemann problem with this initial data
is known as modified Sod's shock tube \cite{toro2009} and is specifically designed to have a sonic point within the rarefaction wave.
Another difference to the classical shock tube problem is that the left boundary is now an inlet.

\begin{figure}[ht!]
\centering
\begin{subfigure}[b]{0.49\textwidth}
\caption{GT-BP}
\includegraphics[width=\textwidth]{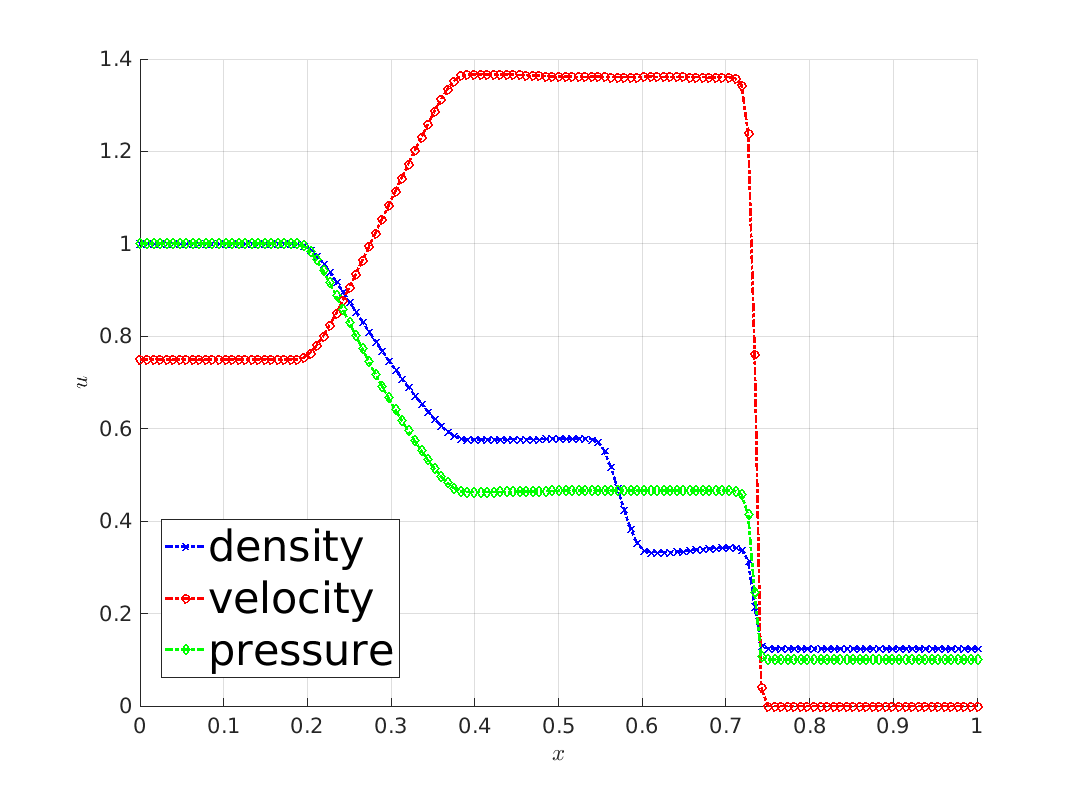}
\end{subfigure}
\begin{subfigure}[b]{0.49\textwidth}
\caption{RT-BP}
\includegraphics[width=\textwidth]{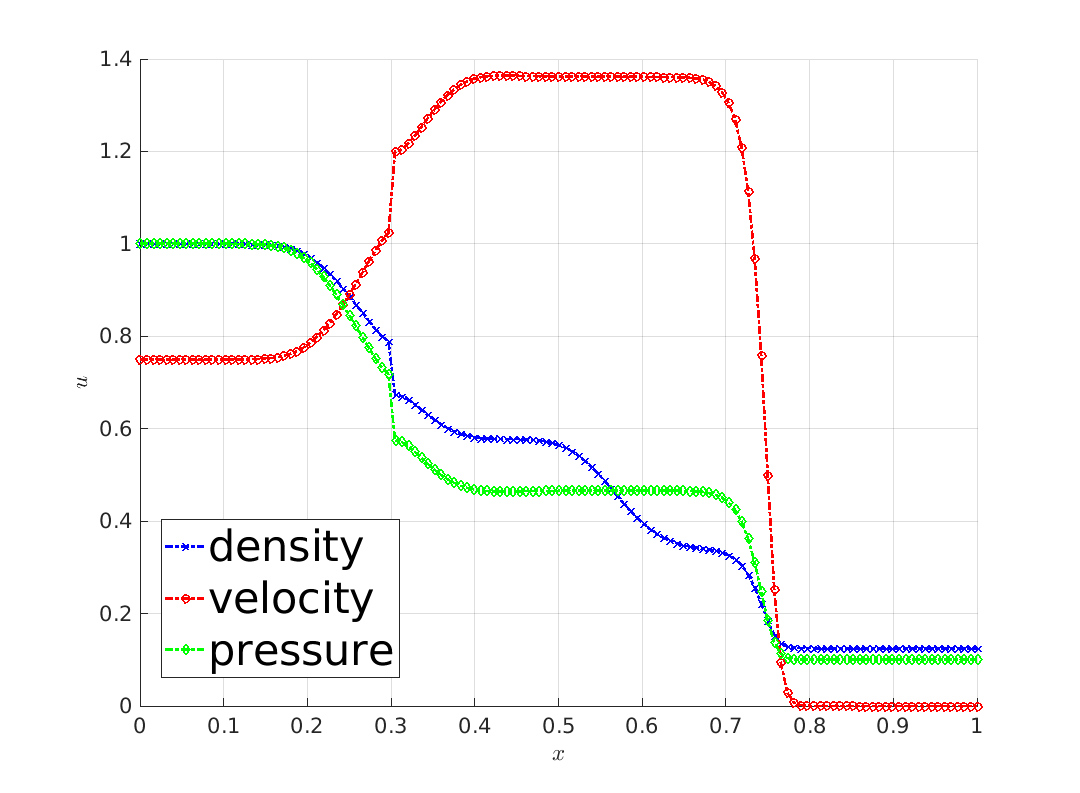}
\end{subfigure}
\begin{subfigure}[b]{0.49\textwidth}
\caption{RT-SD-ED}
\includegraphics[width=\textwidth]{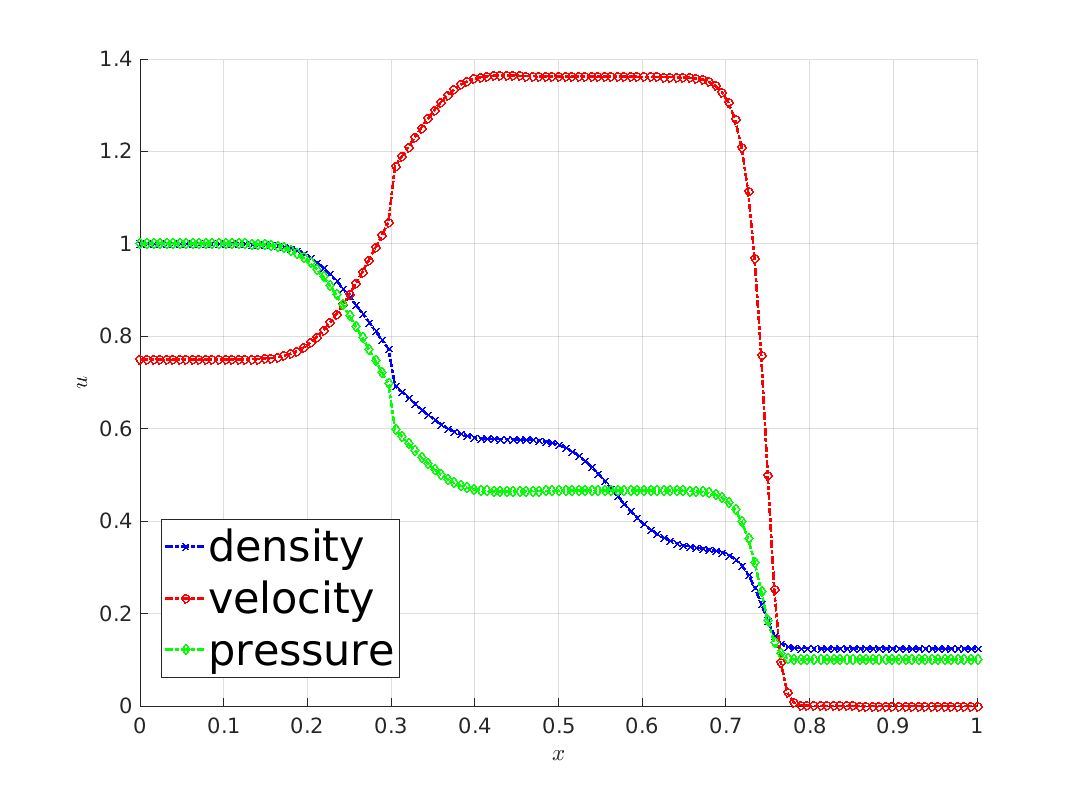}
\end{subfigure}
\begin{subfigure}[b]{0.49\textwidth}
\caption{RT-FDE-ED}
\includegraphics[width=\textwidth]{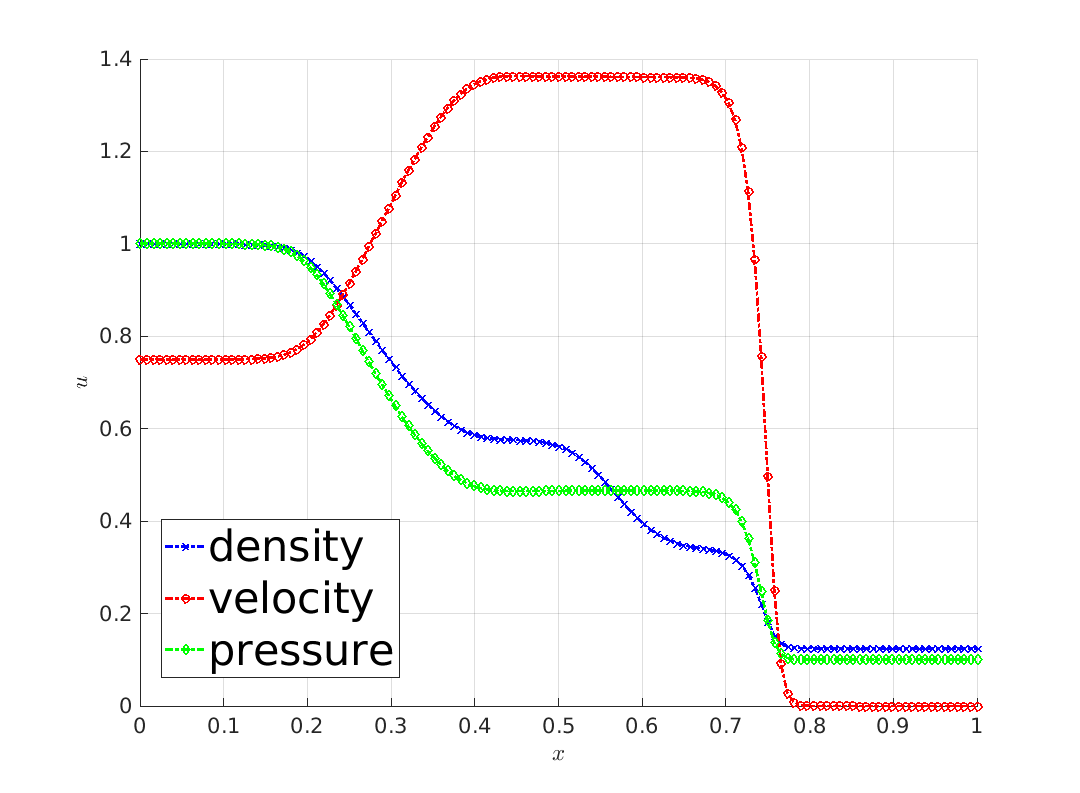}
\end{subfigure}
\caption{Modified Sod's shock tube \cite{toro2009}, primitive variables at $t=0.2$
  calculated on a uniform mesh with $h=\frac 1 {128}$ and $\Delta t= 10^{-3}$.
  No entropy fix is performed 
  in the BP version. The SD and FDE entropy fixes for Roe target fluxes
  use ED bounds \eqref{Qij-ED}.}\label{fig:modsod}
\end{figure}

The results in \cref{fig:modsod} depict approximations to the primitive variables $\rho,v$, and $E$ at $t=0.2$. No entropy fixes are used for Galerkin target fluxes. As in the dam break example, the GT-BP results are qualitatively correct, so there is no need for entropy correction. Although Roe's scheme without an entropy fix is considerably more diffusive, it produces an entropy shock at the sonic point, as in the case of the shallow water equations. The semi-discrete fix reduces the magnitude of the spurious jump, and the resulting approximations do converge to the correct solution as the mesh and time step are refined. However, only the fully discrete fix prevents formation of the entropy shock already on coarse meshes.

\subsubsection{Other benchmarks for the Euler equations}
\label{sec:other}

\begin{figure}[ht!]
\centering
\begin{subfigure}[b]{0.49\textwidth}
\caption{\cite{einfeldt1988}, $E_h=200,~\Delta t = 4\cdot 10^{-4}$}\label{fig:einfeldt}
\includegraphics[width=\textwidth]{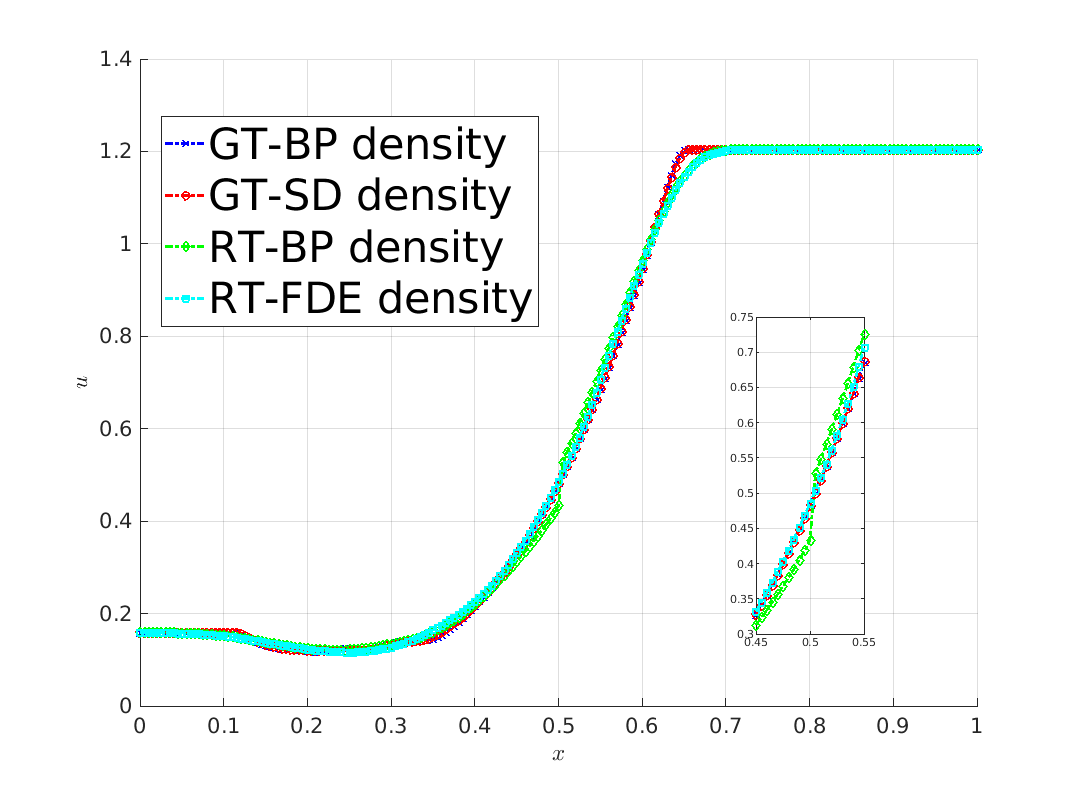}
\end{subfigure}
\begin{subfigure}[b]{0.49\textwidth}
\caption{\cite{moschetta2000}, $E_h=200,~\Delta t = 4\cdot 10^{-4}$}\label{fig:glitch1}
\includegraphics[width=\textwidth]{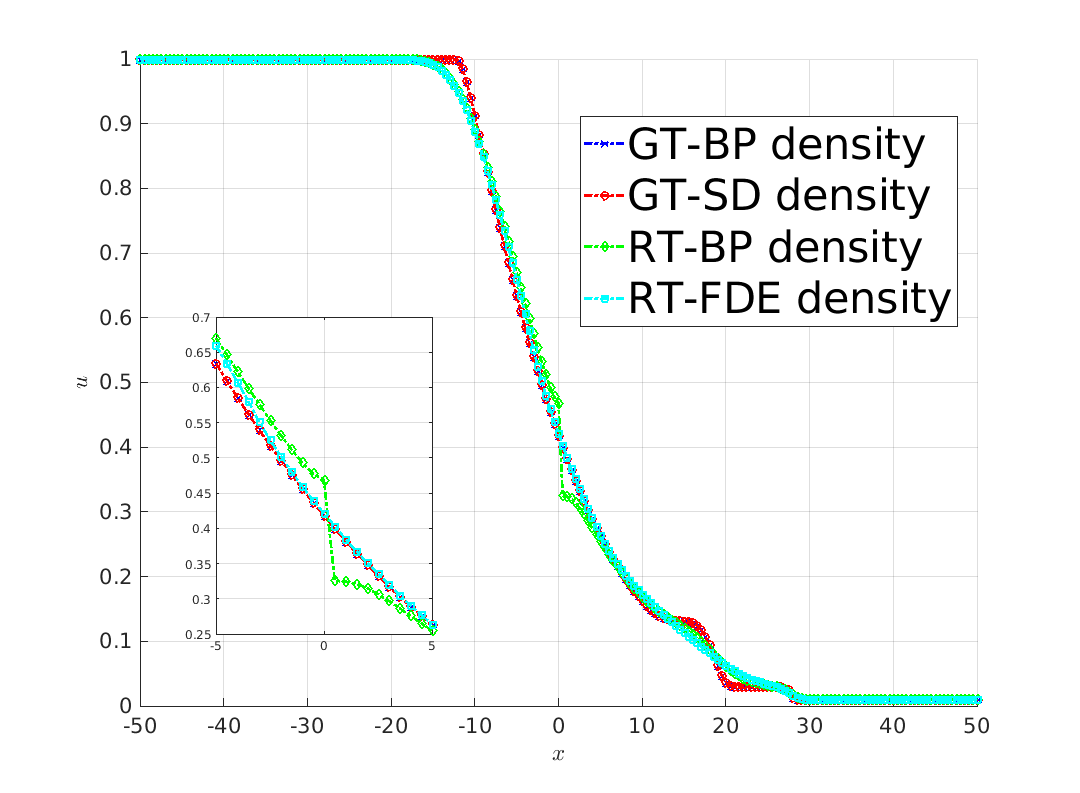}
\end{subfigure}
\begin{subfigure}[b]{0.49\textwidth}
\caption{\cite{moschetta2000}, $E_h=200,~\Delta t = 5\cdot 10^{-6}$}\label{fig:glitch2}
\includegraphics[width=\textwidth]{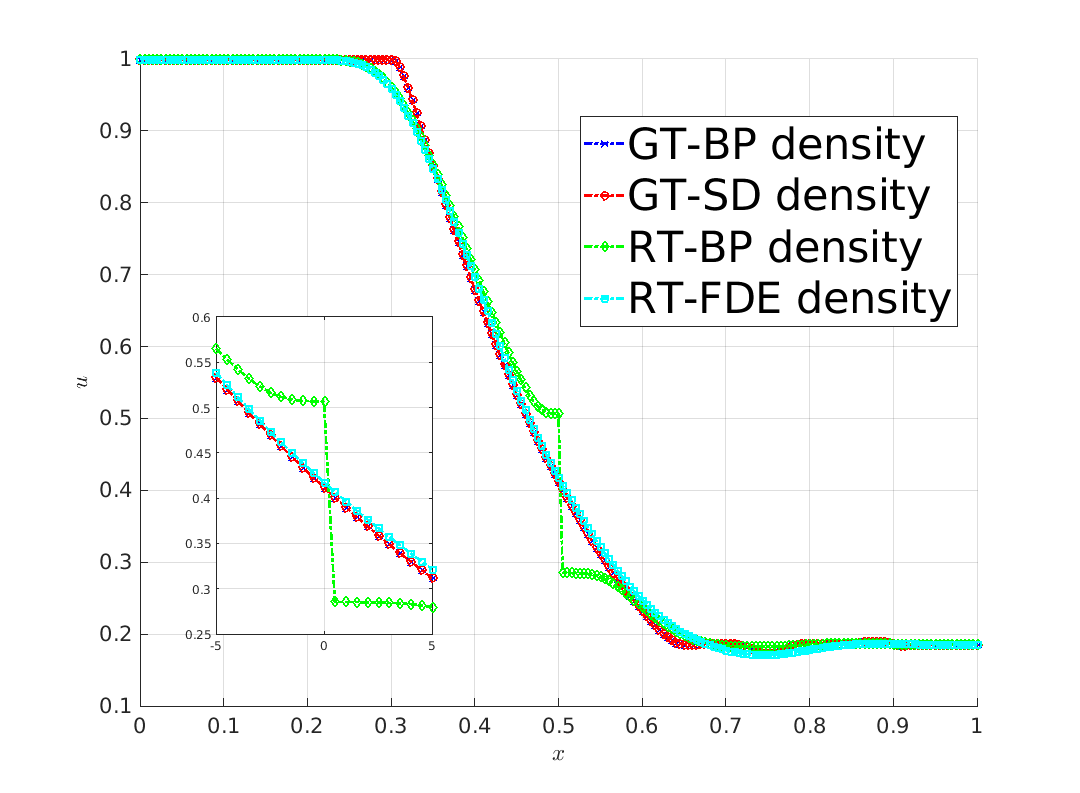}
\end{subfigure}
\begin{subfigure}[b]{0.49\textwidth}
\caption{\cite{shu1989}, $E_h=500,~\Delta t = 2 \cdot 10^{-3}$}\label{fig:sineshock}
\includegraphics[width=\textwidth]{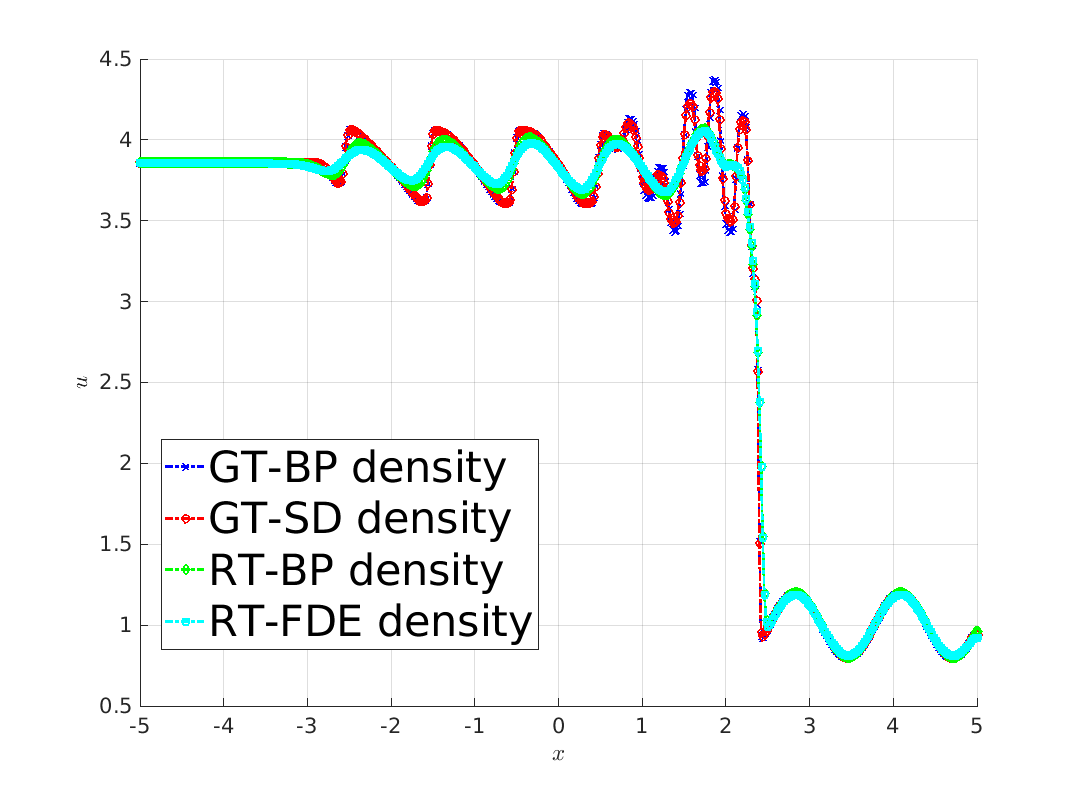}
\end{subfigure}
\caption{Density profiles of the GT-BP, GT-SD, RT-BP, and RT-FDE
  solutions to the 1D Euler equations for standard test problems referenced in the text of Section \ref{sec:other}.}\label{fig:other}
\end{figure}

Among the benchmarks considered so far, only the 1D and 2D KPP problems 
did require entropy fixes for BP schemes using Galerkin fluxes. Indeed, all GT results for hyperbolic systems were qualitatively correct. While the semi-discrete fix did not significantly increase the levels of diffusivity, the fully discrete explicit fix did but was found to be a better cure for entropy shocks generated by Roe target fluxes. To show that this behavior is not unique to the problems considered in \cref{sec:dam} and \cref{sec:sod}, we apply the GT-BP, GT-SD, RT-BP and RT-FDE methods to four additional 1D benchmarks for the Euler equations. The resulting density approximations are presented in \cref{fig:other}. For a detailed description and setup of each test problem, we refer the reader to the references cited below.

The diagram in
\cref{fig:einfeldt} shows the results for the final test in \cite[Sec. 6\,A]{einfeldt1988}.
We remark that the parameter $c_R$ in the setup of this problem is left unspecified in the reference.
Based on the results presented in \cite{einfeldt1988}, we use $c_R=4$.
The results in \cref{fig:glitch1,fig:glitch2} correspond to test cases 1 and 2 in \cite{moschetta2000}. The problem solved in \cref{fig:sineshock} is the classical Shu--Osher sine-shock interaction \cite{shu1989}. Again, approximations obtained with Galerkin target fluxes are free of artifacts, and their accuracy is not affected by the optional semi-discrete fix. The use of Roe target fluxes makes the baseline scheme more diffusive and entropy unstable but the latter deficiency can be cured by applying a limiter-based entropy fix. This study confirms our previous observations that, in practice, the semi-discrete entropy fix for Galerkin fluxes is sufficient to avoid convergence to entropy-violating weak solutions. It does not degrade the overall accuracy and even the entropy-conservative version converges to the vanishing viscosity solution, while reducing the magnitude of spurious jumps on coarse meshes. The fully discrete explicit fix makes it possible to suppress nonphysical effects completely but introduces additional numerical dissipation. Remarkably, the methods under investigation exhibit this behavior for all conservation laws and initial data
considered in this work.

\section{Conclusions}\label{sec:conclusion}

The main outcome of this work is a general framework for
constraining a continuous finite element approximation
to satisfy entropy stability conditions. Combining
a property-preserving
algebraic Lax--Friedrichs method with a high-order
target scheme, we designed algebraic flux correction
procedures that ensure not only preservation of invariant
domains but also validity of local entropy
inequalities. The semi-discrete entropy production
limiter proposed in \cite{kuzmin2020f} was
extended to systems and fully discrete schemes.
The results of our numerical
experiments indicate that a semi-discrete entropy fix is
usually sufficient for convergence to correct weak solutions
if the underlying inequality constraints are formulated using
a dissipative bound for entropy production by antidiffusive
fluxes. The fully discrete explicit fix for forward Euler
stages of an SSP Runge--Kutta method was found to stabilize
an entropy-conservative space discretization
but degrade the rate of convergence to smooth solutions. As
an alternative, we proposed an implicit iterative fix for the
final stage of a general Runge--Kutta method. In terms of accuracy,
this algorithm performs similarly to the semi-discrete fix for
RK stages. However, fully discrete entropy stability is required
by the presented generalization of the Lax--Wendroff theorem to
finite elements. Hence, it is at least a desirable theoretical
property which may need to be enforced if the baseline scheme
is highly unstable. The Galerkin group approximation considered
in this work is almost entropy conservative for $\mathbb{P}_1$
and $\mathbb{Q}_1$ finite elements. For such target schemes, the
benefits of a fully discrete fix may not be worth the effort.
The entropy-conservative version of our semi-discrete fix
has already been extended to high-order continuous finite
element discretizations of scalar conservation laws in
\cite{kuzmin2020g}, while a high-order DG version of the
bound-preserving MCL limiter without entropy fixes was
developed in \cite{hajduk2021}.  We envisage
that entropy stability of such high-order AFC schemes can
also be enforced using the proposed methodology.

\section*{Acknowledgments}

This work was supported by the German Research Association (DFG) under grant KU 1530/23-3.

\appendix

\bibliographystyle{bibstyle-article}
\bibliography{bibliography}
\addcontentsline{toc}{section}{References}
\end{document}